# COEXISTENCE IN LOCALLY REGULATED COMPETING POPULATIONS AND SURVIVAL OF BRANCHING ANNIHILATING RANDOM WALK

By Jochen Blath[1], Alison Etheridge [2] and Mark Meredith[3]

*Technische Universität Berlin, University of Oxford and University of Oxford*

We propose two models of the evolution of a pair of competing populations. Both are lattice based. The first is a compromise between fully spatial models, which do not appear amenable to analytic results, and interacting particle system models, which do not, at present, incorporate all of the competitive strategies that a population might adopt. The second is a simplification of the first, in which competition is only supposed to act within lattice sites and the total population size within each lattice point is a constant. In a special case, this second model is dual to a branching annihilating random walk. For each model, using a comparison with oriented percolation, we show that for certain parameter values, both populations will coexist for all time with positive probability. As a corollary, we deduce survival for all time of branching annihilating random walk for sufficiently large branching rates. We also present a number of conjectures relating to the rôle of space in the survival probabilities for the two populations.

**1. Introduction.** Natural populations interact with one another and with their environment in complex ways. No mathematical model can possibly incorporate all such interactions and still remain analytically tractable. As a result, in order to understand the effects of a feature of a population's dynamics, it is often useful to study "toy models." In this paper, we investigate

Received February 2004; revised March 2007.
[1]Supported by EPSRC GR/R985603.
[2]Supported by EPSRC Advanced Fellowship GR/A90923.
[3]Supported by an EPSRC Doctoral Training Account.
*AMS 2000 subject classifications.* Primary 60K35; secondary 60J80, 60J85, 60J70, 92D25.
*Key words and phrases.* Competing species, coexistence, branching annihilating random walk, interacting diffusions, regulated population, heteromyopia, stepping stone model, survival, Feller diffusion, Wright–Fisher diffusion.

two such toy models that aim to parody the evolution of two populations that are distributed in space and competing for the same resource. Both of our models can be viewed as a compromise between fully spatial models which do not appear to be amenable to a rigorous mathematical analysis and interacting particle system models which do not, at present, incorporate all of the competitive strategies that a population of, say, plants might adopt.

Although lattice based, our first model is highly reminiscent of the models in continuous space studied by Bolker and Pacala [2] and Murrell and Law [12], while admitting a rigorous mathematical analysis. It comprises a system of interacting diffusions, indexed by $\mathbb{Z}^d$, driven by independent Feller noises and coupled through a drift term that reflects migration and competition (both within and between species). Our second model is much simpler: first, we suppose that the parameters governing migration of individuals within the two populations are the same and that competition between the populations acts only within individual lattice sites; second, we suppose that the total population size within each lattice site is a fixed constant. When we further restrict to the *symmetric* case, in which the parameters governing the evolution of the two populations are the same, we shall exhibit a duality between this second process and a *branching annihilating random walk*. The latter is a process that has received considerable attention in the physics literature and we believe this duality to be of some interest in its own right.

A natural starting point for modeling two competing populations is the classical Lotka–Volterra model. This is a deterministic model for the evolution of the total sizes of the two populations, denoted $N_1(t), N_2(t)$. They are assumed to follow the following system of differential equations:

$$(1) \quad \frac{dN_1}{dt} = r_1 N_1 \left(1 - \frac{N_1}{K_1} - \alpha_{12} \frac{N_2}{K_1}\right), \qquad \frac{dN_2}{dt} = r_2 N_2 \left(1 - \frac{N_2}{K_2} - \alpha_{21} \frac{N_1}{K_2}\right),$$

where $r_i$, $K_i$ are respectively the intrinsic growth rates and carrying capacities of the two species and the $\alpha_{ij}$ measure the interspecific competition. It is easy to check that longterm coexistence of the two populations is possible if $K_1 > \alpha_{12} K_2$ and $K_2 > \alpha_{21} K_1$. A number of models have been proposed that extend this in two different ways. First, they incorporate *spatial structure* into the populations and second, they assume that the evolution of the populations is *stochastic*.

It is far from clear how spatial structure affects the chances of longterm coexistence for two competing populations. Traditionally, ecologists have believed that the local nature of interactions between populations that are dispersed in space promotes coexistence. One reason is the so-called *competition–colonization trade-off*: a weaker competitor that is good at colonization may be able to survive by exploiting "gaps" between its competitors. It has also been claimed that because in spatial models the population tends to become segregated into clusters of a single type, the intraspecific competition will be more important than the interspecific competition. Pacala and Levin [14] make an attempt to quantify this effect. On the other hand, Neuhauser and Pacala [13] propose and analyze a spatial stochastic model for competing species in which space actually makes coexistence harder. This suggests, then, that in their model, it is actually the interactions at the cluster boundaries that dominate.

In order to obtain analytic results about spatial stochastic models, simplifying assumptions must be made. Murrell and Law [12] point out that common assumptions are that the parameters of neighborhoods over which individuals compete are the same, irrespective of species, or that dispersal and competition neighborhoods are of the same size, but that dropping such symmetries can have profound consequences. They argue, using a simulation study and the method of moment closure for a specific stochastic model in two space dimensions, that spatial structure *can* promote coexistence, by showing that in the spatial setting, two populations in which the overall strength of interspecific and intraspecific competition is the same can coexist, but only if the distance over which individuals sense their heterospecific neighbors (i.e., their competitors) is shorter than that over which they sense their conspecific neighbors. They coin the term *heteromyopia* for populations that are "shortsighted" in this way. We explain this concept in a little more detail in the context of our first model in Section 2. Although this model admits such differences in neighborhood size, our methods are not strong enough to confirm the numerical findings of Murrell and Law in this context. Indeed, even when the populations migrate in a symmetric way and intraspecific and interspecific competition neighborhoods are of the same size, although we conjecture (in Section 2) that space does *not* make coexistence harder for our model, our methods are not strong enough to provide a rigorous proof of this claim.

*Model* I. Following Bolker and Pacala [2], we assume that the strategies for survival that individuals in our model can employ are (i) to colonize relatively unpopulated areas quickly, (ii) to quickly exploit resources in those areas and (iii) to tolerate local competition. We take two different populations (species) and each can adopt a different combination of strategies for survival. In order to simplify the proofs of our results, we suppose our populations to be living on $\mathbb{Z}^d$ (the biologically relevant case is $d=2$). The dynamics of the model are entirely analogous to those considered by Bolker and Pacala [2] and by Murrell and Law [12]. We write $\{X(t)\}_{t\geq 0} = \{X_i(t), i \in \mathbb{Z}^d\}_{t\geq 0}$ and $\{Y(t)\}_{t\geq 0} = \{Y_i(t), i \in \mathbb{Z}^d\}_{t\geq 0}$ for our two populations. We shall suppose that the pair of processes $\{X(t)\}_{t\geq 0}$, $\{Y(t)\}_{t\geq 0}$ solves the following system of stochastic differential equations:

$$dX_i(t) = \sum_{j \in \mathbb{Z}^d} m_{ij}(X_j(t) - X_i(t))\,dt$$

$$\text{(2)} \qquad + \alpha\left(M - \sum_{j \in \mathbb{Z}^d} \lambda_{ij} X_j(t) - \sum_{j \in \mathbb{Z}^d} \gamma_{ij} Y_j(t)\right) X_i(t)\, dt$$

$$+ \sqrt{\sigma X_i(t)}\, dB_i(t),$$

$$dY_i(t) = \sum_{j \in \mathbb{Z}^d} m'_{ij}(Y_j(t) - Y_i(t))\, dt$$

$$\text{(3)} \qquad + \alpha'\left(M' - \sum_{j \in \mathbb{Z}^d} \lambda'_{ij} Y_j(t) - \sum_{j \in \mathbb{Z}^d} \gamma'_{ij} X_j(t)\right) Y_i(t)\, dt$$

$$+ \sqrt{\sigma Y_i(t)}\, dB'_i(t),$$

where $\{\{B_i(t)\}_{t \geq 0}, \{B'_i(t)\}_{t \geq 0}, i \in \mathbb{Z}^d\}$ is a family of independent standard Brownian motions. The (bounded nonnegative) parameters $m_{ij}$, $m'_{ij}$, $\lambda_{ij}$, $\lambda'_{ij}$, $\gamma_{ij}$ and $\gamma'_{ij}$ are all supposed to be functions of $\|i - j\|$ alone and to vanish for $\|i - j\| > R$ for some $R < \infty$. In other words, the range of both migration and interaction for the two populations will be taken to be finite.

Here, $\|\cdot\|$ can either denote the lattice distance (so that simple random walk is included) or the maximum norm on $\mathbb{Z}^d$, but it will be convenient to take the maximum norm. Moreover, notice that by a change of units, there is no loss of generality in taking the same $\sigma$ for both populations and indeed, hereafter, we may and will set $\sigma = 1$.

REMARK (Existence and uniqueness). Note that Model I is not covered by the now standard results in Shiga and Shimizu [17]. However, Blath, Etheridge and Meredith [1] provide, for $1 \leq p \in \mathbb{N}$, *existence* of a continuous positive (weak) solution in the space $\ell_\Gamma^{4p} = \{x \in \mathbb{R}^{\mathbb{Z}^d} : \|x\|_{\Gamma, 4p} < \infty\}$, where the weighted $\ell^p$-norm $\|\cdot\|_{\Gamma, p}$ is defined, for $\Gamma = \{\Gamma_i\}_{\mathbb{Z}^d} \in l^1((0, \infty)^{\mathbb{Z}^d})$, by $\|x\|_{\Gamma, p} = (\sum_{i \in \mathbb{Z}^d} \Gamma_i |x_i|^p)^{1/p}$. We assume $\Gamma_i / \Gamma_j < f(\|i - j\|)$ for some continuous function $f : [0, \infty) \to [0, \infty)$. For example, set $\Gamma_i = e^{-\|i\|}$ for each $i \in \mathbb{Z}^d$. *Uniqueness* remains open, after considerable efforts, including those of several experts whom we have consulted. At first sight, one expects to be able to prove uniqueness in a suitable weighted $\ell_\Gamma^{4p}$ space by an application of the (infinite-dimensional) Yamada–Watanabe theorem. This works only in the special case when $\lambda_{ij}$ and $\gamma_{ij}$ both vanish for $i \neq j$. The nonlocal nature of the interaction destroys the vestiges of monotonicity available in this special case.

REMARK (Blath, Etheridge and Meredith [1]). The full version of this paper, Blath, Etherige and Meredith [1], which has also successfully undergone the peer reviewing process of *Annals of Applied Probability*, proved too long to be published in its entirety. It contains full technical details and some additional remarks and is available from the webpages of the authors.

For the $X$-population, the first two strategies for survival listed above correspond to taking large $m_{ij}$ and large $\alpha M$, while the third corresponds to taking small $\lambda_{ij}$ (conspecific competition) and $\gamma_{ij}$ (interspecific competition). By varying $M$, we can also model how efficiently the species uses the available resources: a species that can tolerate lower resource levels will have a higher value of $M$.

DEFINITION 1.1 (Notions of survival). Let $p \in [0,1)$. We shall say that the $X$-population *survives for all time* with probability greater than $p$ if there exists $\kappa > 0$ such that
$$\liminf_{t \to \infty} \mathbb{P}[X_0(t) > \kappa] > p.$$
We shall say that *both populations persist for all time* with probability greater than $p$ if there exists $\kappa > 0$ such that
$$\mathbb{P}[\forall t > 0, \exists i, j \in \mathbb{Z}^d : X_i(t), Y_j(t) > \kappa] > p.$$
Finally, we shall say that the populations exhibit *long-term coexistence* with probability greater than $p$ if there exists $\kappa > 0$ such that
$$\liminf_{t \to \infty} \mathbb{P}[X_0(t), Y_0(t) > \kappa] > p.$$

Observe that the third notion is much stronger than the second one. Also, note that if $\gamma_{ij} = \gamma'_{ij}$ is zero for all $i, j \in \mathbb{Z}^d$, then each population follows an independent copy of the so-called *stepping stone version of the Bolker–Pacala model* introduced in Etheridge [9]. There, it is proved that if the range of migration is at least as great as the range over which the population interacts with itself (here determined by the $\{\lambda_{ij}\}$), then provided that $\alpha M$ is sufficiently large, the population will survive. A partial converse of this, proved there only in the context of a continuous-space analogue of this model, suggests that this condition is actually necessary, a conclusion reached independently by Law, Murrell and Dieckmann [11]. We therefore assume from the outset that there exists a constant $c > 0$ such that for all $i, j \in \mathbb{Z}^d$, we have $m_{ij} > c\lambda_{ij}$ (resp. $m'_{ij} > c\lambda'_{ij}$) whenever $\lambda_{ij}$ (resp. $\lambda'_{ij}$) is strictly positive. Indeed, Theorem 1.5 in Etheridge [9] then tells us that if $\alpha M > \sum_j m_{ij}$ is sufficiently large (depending on $c$), the single species model for $X$ started from any nontrivial translation invariant initial condition survives with positive probability, that is, there exists a $\kappa > 0$ such that $\liminf_{t \to \infty} \mathbb{P}[X_0(t) > \kappa] > 0$.

For the competing species model, we will have to make similar and additional assumptions. In particular, we shall choose initial conditions in such a way that we can find a box where both populations are present, but not so prevalent that the competitive interaction between them is too large.

*Notation and assumptions for Theorem* 1.2.

- The parameters $m_{ij}$, $m'_{ij}$, $\lambda_{ij}$, $\lambda'_{ij}$, $\gamma_{ij}$ and $\gamma'_{ij}$ are nonnegative functions of $\|i-j\|$ alone and vanish for $\|i-j\| > R$ for some $R < \infty$.
- $\{m_{ij}\}$, $\{m'_{ij}\}$, $\{\lambda_{ij}\}$ and $\{\lambda'_{ij}\}$ are fixed in such a way that there exists a constant $c > 0$ such that, for all $m_{ij}, m'_{ij} \neq 0$,

(4) $$\frac{1}{c}\lambda_{ij} < m_{ij} < c\lambda_{ij} \quad \text{and} \quad \frac{1}{c}\lambda'_{ij} < m'_{ij} < c\lambda'_{ij}.$$

For all $i,j$ such that $m_{ij} = 0$ (resp. $m'_{ij} = 0$) we require $\lambda_{ij} = 0$ (resp. $\lambda'_{ij} = 0$). Assume that $\{m_{ij}\}$ and $\{m'_{ij}\}$ are nondiagonal and of the same range and that $\lambda_{ii}, \lambda'_{ii} > 0$ for all $i \in \mathbb{Z}^d$.

- Let $L = \max\{\|j - i\| : m_{ij}, m'_{ij} \neq 0\} \le R < \infty$.
- Assume that $\alpha M > \sum_{j \in \mathbb{Z}^d} m_{ij}$ and $\alpha' M' > \sum_{j \in \mathbb{Z}^d} m'_{ij}$.
- Let $b \in \mathbb{Z}$ such that $\max\{\|j - i\| : \gamma_{ij} \text{ or } \gamma'_{ij} \neq 0\}$ is less than $(b-1)L$.
- For $m \in \mathbb{N} \cup \{\infty\}$ and $0 < \kappa_1 < \kappa_2 < \infty$, $0 < \kappa'_1 < \kappa'_2 < \infty$, we write

$$(X(0), Y(0)) \in H(\kappa_1, \kappa_2; \kappa'_1, \kappa'_2; m),$$

if $X(0), Y(0) \in \ell_\Gamma^{4p}$ and there exists a box $J = \{[-m,m]^d \cap \mathbb{Z}^d\} \subset \mathbb{Z}^d$, such that for all $i \in J$,

(5) $$X_0(i) \in [\kappa_1, \kappa_2) \quad \text{and} \quad Y_0(i) \in [\kappa'_1, \kappa'_2).$$

REMARK. One can drop the assumption that the range of $\{m_{ij}\}$ and $\{m'_{ij}\}$ are the same, but this will make the proof much more awkward.

THEOREM 1.2. *Under the above assumptions, there exist finite constants $M_0 > 0, M'_0 > 0$ such that:*

1. *for each $M > M_0$ and $M' > 0$, there is a constant $\gamma = \gamma(M, M') > 0$ and constants $0 < \kappa_1 < \infty$, $0 < \kappa'_2 < \infty$, such that if $\sum_j \gamma_{ij} < \gamma$ and*

   $$(X(0), Y(0)) \in H(\kappa_1, \infty; 0, \kappa'_2; (b+1/2)L),$$

   *then the $X$-population survives for all time with probability greater than one half;*

2. *similarly, for each $M' > M'_0$ and $M > 0$, there is a constant $\gamma' = \gamma'(M, M') > 0$ and constants $0 < \kappa_2 < \infty$, $0 < \kappa'_1 < \infty$, such that if $\sum_j \gamma'_{ij} < \gamma'$ and*

   $$(X(0), Y(0)) \in H(0, \kappa_2; \kappa'_1, \infty; (b+1/2)L),$$

   *then the $Y$-population survives for all time with probability greater than one half.*

COROLLARY 1.3. *Under the conditions of Theorem 1.2, for each pair $(M, M')$ with $M > M_0$ and $M' > M'_0$, there is a pair $(\gamma, \gamma')$ with $\gamma > 0, \gamma' > 0$ and constants $0 < \kappa_1 < \kappa_2 < \infty$, $0 < \kappa'_1 < \kappa'_2 < \infty$, such that if $\sum_j \gamma_{ij} < \gamma$, $\sum_j \gamma'_{ij} < \gamma'$ and*

$$(X(0), Y(0)) \in H(\kappa_1, \kappa_2; \ \kappa'_1, \kappa'_2; \ (b+1/2)L),$$

*then the $X$- and $Y$-populations both persist for all time with positive probability. Moreover, for each such pair, if $H(\kappa_1, \kappa_2; \kappa'_1, \kappa'_2; (b+1/2)L)$ is replaced by $H(\kappa_1, \kappa_2; \kappa'_1, \kappa'_2; \infty)$, then there is longterm coexistence with positive probability, that is, there exists $\kappa > 0$ such that*

$$\liminf_{t \to \infty} \mathbb{P}[X_0(t), Y_0(t) > \kappa] > 0.$$

As we will explain in Section 2, we would conjecture a result very much stronger than Theorem 1.2 (or Corollary 1.3). In particular, we provide evidence to support the claim that in the biologically relevant case of two dimensions, if we take the special case of our model in which $\alpha = \alpha'$, $M = M'$ and $m_{ij} = m'_{ij}$, then provided that $\gamma'_{ij} \leq \lambda_{ij}$ (resp. $\gamma_{ij} \leq \lambda'_{ij}$) with strict inequality whenever $\lambda_{ij}$ (resp. $\lambda'_{ij}$) $\neq 0$, and that the parameters are such that if $\gamma_{ij}$ and $\gamma'_{ij}$ were zero then the single species models would survive, then with positive probability, the competing species model will coexist for all time. This would be precisely the prediction of the corresponding Lotka–Volterra model. If we drop the assumptions $\alpha = \alpha'$ and $M = M'$, then this conjecture must be modified to reflect competition–colonization trade-off. We formulate this and other conjectures more carefully in Section 2. In the process, we are led to consider our second model of two competing species.

*Model* II. Suppose now that the neighborhood over which each individual competes is just the site in which it lives, so that the only interaction between different points in $\mathbb{Z}^d$ is through migration. In addition, we suppose that the migration mechanism for the two populations is the same and that the *total* population size in each site is constant [i.e., $X_i(t) + Y_i(t) \equiv N > 0$ for all $i \in \mathbb{Z}^d$ and all $t \geq 0$]. Let us write $p_i(t) = X_i(t)/N$ for the proportion of the total population in $i$ at time $t$ that belongs to the $X$-population. Then, as we will see in Section 2, we arrive at the much simpler model

$$\begin{aligned}dp_i(t) = \sum_{j \in \mathbb{Z}^d} m_{ij}(p_j(t) - p_i(t))\, dt &+ sp_i(t)(1-p_i(t))(1-\mu p_i(t))\, dt \\ &+ \sqrt{N^{-1} p_i(t)(1-p_i(t))}\, dW_i(t),\end{aligned} \quad (6)$$

where $s = \alpha M - \alpha' M' + (\alpha' \lambda'_{ii} - \alpha \gamma_{ii})N$ and

$$\mu = \frac{(\alpha' \lambda'_{ii} - \alpha \gamma_{ii})N + (\alpha \lambda_{ii} - \alpha' \gamma'_{ii})N}{\alpha M - \alpha' M' + (\alpha' \lambda'_{ii} - \alpha \gamma_{ii})N},$$

and, finally, $\{W_i(t), i \in \mathbb{Z}^d\}_{t \geq 0}$ is a family of independent Brownian motions. This model is a system of interacting Fisher–Wright diffusions for gene frequencies in a spatially structured population. From the results in Shiga and Shimizu [17], it follows that if $p_i(0) \in [0,1]$ for all $i \in \mathbb{Z}^d$, then this system has a continuous, pathwise unique, $[0,1]^{\mathbb{Z}^d}$-valued strong solution for all times $t \geq 0$.

If $\mu < 1$, then in each site $i$, there is selection in favor of either the $X$-type or the $Y$-type according to whether $s > 0$ or $s < 0$. If $\mu > 1$, then in each site $i$, we have selection in favor of *heterozygosity* if $s > 0$ and selection in favor of *homozygosity* if $s < 0$. In the "neutral" case ($s = 0$), the process has a *moment dual*, the so-called *structured coalescent* (see, e.g., Shiga [16]), and it is easy to show that if $d \geq 3$, then with positive probability, there will be longterm coexistence of our two populations, whereas if $d \leq 2$, with probability one, eventually only one population will be present.

Notice that we have selection in favor of heterozygosity precisely when

$$(\alpha \lambda_{ii} - \alpha' \gamma'_{ii})N > \alpha M - \alpha' M' \quad \text{and} \quad (\alpha' \lambda'_{ii} - \alpha \gamma_{ii})N > \alpha' M' - \alpha M.$$

We sketch a proof of the following result and present a more detailed analysis in a forthcoming work.

THEOREM 1.4. *Let $\{p_i(t), i \in \mathbb{Z}^d\}_{t \geq 0}$ evolve according to Model* II. *Suppose that $\mu > 1$ and let $\varepsilon \in (0, 1/4)$. Then, if $p_i(0) \in (\varepsilon, 1-\varepsilon)$ for some $i \in \mathbb{Z}^d$, there exists an $s_0 \in [0, \infty)$ such that for all $s > s_0$, we have*

$$\mathbb{P}[\forall t > 0, \exists i \in \mathbb{Z}^d : \varepsilon < p_i(t) < 1 - \varepsilon] > 0.$$

*Moreover, if $p_i(0) \in (\varepsilon, 1-\varepsilon)$ for all $i \in \mathbb{Z}^d$, then*

$$\liminf_{t \to \infty} \mathbb{P}[\varepsilon < p_0(t) < 1 - \varepsilon] > 0.$$

In the case when the two populations evolve symmetrically, that is, $\mu = 2$, Model II reduces to

$$\begin{aligned}dp_i(t) = \sum_j m_{ij}(p_j(t) - p_i(t))\,dt + sp_i(t)(1-p_i(t))(1-2p_i(t))\,dt \\ + \sqrt{N^{-1}p_i(t)(1-p_i(t))}\,dW_i(t).\end{aligned} \quad (7)$$

For general $s$, there is no convenient dual, but in Lemma 2.1, we find an alternative duality with a system of *branching annihilating random walks*.

DEFINITION 1.5 (Branching annihilating random walk). The Markov process $\{n_i(t), i \in \mathbb{Z}^d\}_{t \geq 0}$ with values $n_i(t) \in \mathbb{Z}_+$ and dynamics described by

$$\begin{cases} n_i \mapsto n_i - 1 \\ n_j \mapsto n_j + 1 \end{cases} \quad \text{at rate } n_i m_{ij} \quad \text{(migration)},$$

COEXISTENCE IN COMPETING POPULATIONS AND BARW        9

$$n_i \mapsto n_i + 2 \quad \text{at rate } sn_i \quad \text{(branching)},$$

$$n_i \mapsto n_i - 2 \quad \text{at rate } \tfrac{1}{2}n_i(n_i - 1) \quad \text{(annihilation)}$$

is called a *branching annihilating random walk* with *branching rate s* (and offspring number two).

COROLLARY 1.6. *There exists $s_0 \geq 0$ such that if $s > s_0$, the branching annihilating random walk, started from an even number of particles at time zero, will survive for all time with positive probability.*

REMARK. Note that in our branching annihilating random walk, a birth event results in one individual splitting into *three*, a net increase of two, whereas an annihilation event results in the loss of two particles. As a result, we have parity preservation: if we start from an odd number of particles, then there will always be an odd number of particles in the system (so, in particular, at least one). This is why we restrict the initial number of particles in Corollary 1.6 to be even.

Branching annihilating random walk has received considerable attention from physicists (see Täuber [18] for a review). For example, Cardy and Täuber [3, 4] consider precisely the process described above. Our conjecture for Model II, stated in Section 2, is based on their results, which, in turn, are based on perturbation theory and renormalization group calculations.

The rest of the paper is organized as follows. In Section 2, we explain the relationship between our two models and the duality between the symmetric form of Model II and branching annihilating random walk. We also make some conjectures about the longterm behavior of our two models and relate them to results and conjectures for other toy models. The proof of our main result will rely upon a comparison with oriented $2N$-dependent percolation and so in Section 3, we recall the definition of $2N$-dependent percolation and state a suitable comparison result. The proofs of Theorem 1.2 and Corollary 1.3 are in Section 4 and a sketch of the proof of Theorem 1.4 is in Section 4.2.3. Corollary 1.6 will then be an immediate consequence of the duality of Model II and branching annihilating random walk.

## 2. Heuristics, duality and relation to existing models.

2.1. *Relationship between the two models.* Suppose that the evolution of our population follows Model I, that is, it is determined by equations (2) and (3). We now derive the system of equations governing the proportion of the total population at time $t$ at site $i$ that belongs to the $X$-subpopulation. We need some notation. If we write $N_i(t) = X_i(t) + Y_i(t)$ and $p_i(t) = X_i(t)/N_i(t)$, then an application of Itô's formula (and some rearrangement) gives

$$dp_i(t) = \sum_{j \in \mathbb{Z}^d} m_{ij} \frac{N_j(t)}{N_i(t)} (p_j(t) - p_i(t)) \, dt$$

$$- \sum_{j \in \mathbb{Z}^d} (m_{ij} - m'_{ij}) p_i(t)(1 - p_i(t)) \, dt$$

$$+ \sum_{j \in \mathbb{Z}^d} (m_{ij} - m'_{ij}) \frac{N_j(t)}{N_i(t)} p_i(t)(1 - p_j(t)) \, dt$$

(8)
$$+ \left[ \alpha M - \alpha' M' + \sum_{j \in \mathbb{Z}^d} (\alpha' \lambda'_{ij} - \alpha \gamma_{ij}) N_j(t) \right.$$

$$\left. + \sum_{j \in \mathbb{Z}^d} (\alpha \gamma_{ij} + \alpha' \gamma'_{ij} - \alpha' \lambda'_{ij} - \alpha \lambda_{ij}) N_j(t) p_j(t) \right]$$

$$\times p_i(t)(1 - p_i(t)) \, dt$$

$$+ \sqrt{N_i(t)^{-1} p_i(t)(1 - p_i(t))} \, dW_i(t),$$

where $\{W_i(t), i \in \mathbb{Z}^d\}_{t \geq 0}$ is a family of independent Brownian motions.

We concentrate on the case when $m_{ij} = m'_{ij}$. Notice that if we also assume that $\lambda_{ij}, \lambda'_{ij}$ and $\gamma_{ij}, \gamma'_{ij}$ are zero for $i \neq j$ and that the population sizes $N_i(t)$ are, in fact, a fixed constant, then we arrive at Model II:

$$dp_i(t) = \sum_j m_{ij}(p_j(t) - p_i(t)) \, dt + s p_i(t)(1 - p_i(t))(1 - \mu p_i(t)) \, dt$$

$$+ \sqrt{N^{-1} p_i(t)(1 - p_i(t))} \, dW_i(t),$$

where $s = \alpha M - \alpha' M' + N(\alpha' \lambda'_{ii} - \alpha \gamma_{ii})$ and

$$\mu = \frac{(\alpha' \lambda'_{ii} - \alpha \gamma_{ii}) N + (\alpha \lambda_{ii} - \alpha' \gamma'_{ii}) N}{\alpha M - \alpha' M' + (\alpha' \lambda'_{ii} - \alpha \gamma_{ii}) N}.$$

2.2. *Conjectures for Model* II. Our conjectures for Model II are based on the symmetric case, implying that $\mu = 2$. The model then reduces to the system (7). In this case, we are able to find a convenient dual process. First, we transform the equations. Let $x_i(t) = 1 - 2p_i(t)$. Then

(9)
$$dx_i(t) = \sum_j m_{ij}(x_j(t) - x_i(t)) \, dt$$

$$+ \frac{s}{2}(x_i^3(t) - x_i(t)) \, dt - \sqrt{(1 - x_i^2(t))} \, dW_i(t).$$

LEMMA 2.1. *The system* (9) *is dual to branching annihilating random walk with branching rate* $s/2$, *denoted* $\{n_i(t), i \in \mathbb{Z}^d\}_{t \geq 0}$, *through the duality relationship*

$$\mathbb{E}[\underline{x}(t)^{\underline{n}(0)}] = \mathbb{E}[\underline{x}(0)^{\underline{n}(t)}],$$

*where* $\underline{x}^{\underline{n}} \equiv \prod_{i \in \mathbb{Z}^d} x_i^{n_i}$.

The proof is completely standard (see e.g. Shiga [15]) and is therefore omitted.

Cardy and Täuber [3, 4] consider the branching annihilating random walk model of Definition 1.5. In particular, their results suggest that although in one dimension the optimal value for $s_0$ in Corollary 1.6 is strictly positive, in two dimensions one can take $s_0 = 0$.

CONJECTURE 2.2. *For Model* II *with* $\mu = 2$ *and* $d = 1$, *there exists a critical value* $s_0 > 0$ *such that the populations described by system* (7) *will both persist for all time with positive probability if and only if* $s > s_0$. *In* $d = 2$, *there is positive probability that both populations will persist for all time if and only if* $s > 0$. *For* $d \geq 3$, *this probability is positive if and only if* $s \geq 0$.

Roughly speaking, for $d \geq 2$, if there is a homozygous advantage, then the population will initially form homogenic clusters, but ultimately it will be the interactions at the cluster boundaries that dominate and one type will go extinct. In the heterozygous advantage case, there will be long term coexistence of species. In one dimension, the heterozygous advantage must be "sufficiently strong" if we are to see coexistence.

In fact, we would go further. In view of the genetic interpretation of Model II, it would be odd if the case $\mu = 2$ were pathological, so we expect that in $d \geq 2$, we will have coexistence for any $s > 0$, $\mu > 1$.

CONJECTURE 2.3. *Conjecture* 2.2 *holds true for any* $\mu > 1$, *where in one dimension*, $s_0$ *will now also depend on* $\mu$.

If this conjecture is true, then in dimensions greater than one, for

$$(\alpha'\lambda'_{ii} - \alpha\gamma_{ii})N > \alpha'M' - \alpha M \quad \text{and} \quad (\alpha\lambda_{ii} - \alpha'\gamma'_{ii})N > \alpha M - \alpha'M',$$

we have positive probability that both populations survive. Comparing the quantities $\alpha'\lambda'_{ii} - \alpha\gamma_{ii}$ and $\alpha\lambda_{ii} - \alpha'\gamma'_{ii}$ tells us about the relative effectiveness of the $X$- and $Y$-populations as competitors. If the first is smaller, then the $X$-population is a less effective competitor. However, provided that $\alpha M > \alpha'M'$, we can even allow it to be negative and have positive probability of survival for the $X$-population. This reflects a competition–colonization trade-off.

2.3. *Conjectures for Model* I. We now turn to Model I. We assume that the migration mechanisms governing the two populations are the same. Suppose first that $\alpha = \alpha'$, $M = M'$, $\lambda_{ij} = \gamma'_{ij}$ and $\lambda'_{ij} = \gamma_{ij}$. We then see that the system of equations (8) looks like a selectively neutral stepping stone model with variable population sizes in each lattice site. If we condition on the trajectories of those population sizes, then this process will have a dual process: a system of coalescing random walks in a space-and-time-varying environment. Showing that there is no longterm coexistence of types amounts to showing that two independent random walks evolving in this environment will, with probability one, eventually meet and coalesce. If the environment is sufficiently well behaved, then one might expect this to be true. Problems will arise if the environment develops large "holes," so that the walkers never meet, or very dense clumps, so that when the walkers do meet, they do so in such a heavily populated site that they do not coalesce before moving apart again. Much of our proof of Theorem 1.2 is devoted to showing that the environment does not clump and a special case of that result says that provided both populations are initially present in sufficient numbers in all sites, the probability that any given site is in a "hole" at time $t$ is uniformly bounded below. We therefore conjecture that in the neutral case, Model I will behave qualitatively in the same way as Model II. In the biologically relevant case of two spatial dimensions, we have been unable to produce a proof.

More generally, we believe, still assuming that $m_{ij} = m'_{ij}$ and $\alpha = \alpha'$, $M = M'$, provided that at least one population persists, the question of longterm coexistence of the populations described by Model I will not be changed by assuming that competition only acts within individual lattice sites and, moreover, in that case, the question of coexistence will be the same as for the populations described by Model II. Namely, we make the following conjecture.

CONJECTURE 2.4. *Let $m_{ij} = m'_{ij}$, $\alpha = \alpha'$, $M = M'$ be fixed. Suppose that both $X$- and $Y$-populations start from nontrivial translation invariant initial conditions and that the parameters are such that each population has positive chance of survival in the absence of the other. Further, let $\lambda_{ij} = \lambda'_{ij}, \gamma_{ij} = \gamma'_{ij}$.*

1. *If $\lambda_{ij} < \gamma_{ij}$ for all $j$, then eventually, only one population will be present.*
2. *If $\lambda_{ij} > \gamma_{ij}$ for all $j$, then if $d \geq 2$, with positive probability, the populations will exhibit longterm coexistence. In one dimension, the same result will hold true provided that $\lambda_{ij} - \gamma_{ij}$ is sufficiently large.*
3. *If $\lambda_{ij} = \gamma_{ij}$ and $d \geq 3$, then with positive probability, both populations will exhibit longterm coexistence.*

   *If $d \leq 2$, then with probability one, one of the populations will eventually die out.*

When $\alpha M \neq \alpha' M'$, we would once again expect to see a competition–colonization trade-off.

2.4. *Heteromyopia.* In view of equation (8), it is easy to see that Murrell and Law's heteromyopia might lead to coexistence. They work in a continuous space with the strength of competition between individuals decaying with their distance apart, according to a Gaussian kernel. The analogue of their model in our setting is the symmetric version of Model I with $\lambda_{ij} = \lambda(\|i - j\|)$, $\gamma_{ij} = \gamma(\|i - j\|)$, where the functions $\lambda$ and $\gamma$ are monotone decreasing and $\sum_j \lambda_{ij} = \sum_j \gamma_{ij}$, but the range of $\lambda_{ij}$ is greater than that of $\gamma_{ij}$. We can think of the effect of this as follows. Over small scales, we have homozygous advantage, over larger scales, heterozygous advantage. Again, we expect to see the population forming homogenic clusters, but now the cluster boundaries will be maintained because the heterogeneity there confers an advantage to individuals within the clusters which counteracts the disadvantage to the individuals actually on the boundary. Reversing the sign to give populations with "heterohyperopia" produces the opposite effect. This is not stable and Murrell and Law observe *founder control* in this case, which means that the outcome of the competition is entirely determined by the initial conditions.

CONJECTURE 2.5. *Assume nontrivial translation invariant initial conditions. Suppose that the parameters of Model* I *are symmetric and such that in the absence of the competitor, each population survives for all time with positive probability. For $d \geq 2$, if the populations are heteromyopic, then we will see longterm coexistence with positive probability, whereas in $d = 1$, the populations must be strongly heteromyopic for there to be longterm coexistence with positive probability.*

2.5. *Relation to existing models.*

*The Murrell–Law model.* Our conjectures for Model I are entirely in agreement with the numerical results of Murrell and Law [12]. They analyze a stochastic version of a continuous-space Lotka–Volterra system, similar to ours. The evolution is characterized in terms of moment equations. These moment equations were derived from a stochastic individual-based model by Dieckmann and Law [7]. Although the assumption of a spatially continuous environment is clearly desirable, the price that they pay is that there are very few analytic tools available for the study of the resulting population models, so they use moment closure, assuming in this case a "power-1" closure. In particular, they ignore dynamics of all *spatial* moments beyond order two. In view of the clustering behavior that is characteristic of populations evolving according to spatial branching models in two dimensions, this method has potential pitfalls. In fact, the control of the clumping of the populations that forms an essential part of our proof of Theorem 1.2 also adds considerable credibility to the moment closure technique for these models and hence to the numerical predictions of Murrell and Law.

Neuhauser and Pacala [2] also consider an explicitly spatial stochastic version of the Lotka–Volterra model. Their model is lattice based, but, in contrast to ours, allows only a single individual to live at each lattice site. Moreover, there is instant recolonization so that there will always be exactly one individual at each site in $\mathbb{Z}^d$. This fixed population size makes it more analogous to Model II than to Model I.

DEFINITION 2.6 (Neuhauser–Pacala model). The Markov process $\{\eta_i(t), i \in \mathbb{Z}^d\}_{t \geq 0}$ in which $\eta_i(t) \in \{1, 2\}$ and with dynamics:

1. if $\eta_i(t) = 1$, it becomes 2 at rate $\frac{\lambda f_2}{\lambda f_2 + f_1}(f_1 + \alpha_{12} f_2)$,
2. if $\eta_i(t) = 2$, it becomes 1 at rate $\frac{f_1}{\lambda f_2 + f_1}(f_2 + \alpha_{21} f_1)$,

where $f_k(i) = \frac{|\{j : \eta_j(t) = k : j \in \mathcal{N}_i\}|}{|\mathcal{N}_i|}$ and $\mathcal{N}_i = i + \{j : 0 < \|j\| \leq R\}$, will be said to follow the *Neuhauser–Pacala* (stochastic spatial Lotka–Volterra) model.

The idea is that an individual of type $k$ will die at a rate determined by the proportion of its neighbors that are conspecific plus some constant multiple of the proportion of heterospecific neighbors. Thus, for example, if, in Model I, we took $\lambda_{ij}$ and $\gamma_{ij}$ to have the same range and to be constant on that range, then a small value of $\alpha_{12}$ would correspond to the ratio $\gamma_{ij}/\lambda_{ij}$ being small. The dead individual is immediately replaced by an offspring of one of its neighbors chosen according to a weight that reflects the relative fecundity of the two types. Thus, for example, $\lambda > 1$ would reflect type 2 being more fecund than type 1. In Model I, this would be modeled by taking $\alpha' M' > \alpha M$. Let us recall some results for this model.

THEOREM 2.7 (Neuhauser and Pacala [13], Theorem 1). *Suppose that $\lambda = 1$, $d = 1$ or $2$ and $\alpha_{12} = \alpha_{21} = \alpha$.*

1. *If $\alpha = 0$, then, except for the one-dimensional nearest-neighbor case, product measure with density $1/2$ is the limiting distribution starting from any nontrivial initial distribution.*
2. *If $\alpha$ is sufficiently small (depending on $R$), then coexistence is possible except for the one-dimensional nearest-neighbor case.*

For Model I, a result entirely analogous to part (2) is a special case of Theorem 1.2. For Model II, the analogue is Theorem 1.4. If we believe Conjecture 2.2, then although in $d = 1$ we require the condition "$\alpha$ sufficiently small," in $d = 2$, the corresponding result is true for all $\alpha < 1$. This corresponds to Conjecture 1 of Neuhauser and Pacala [13].

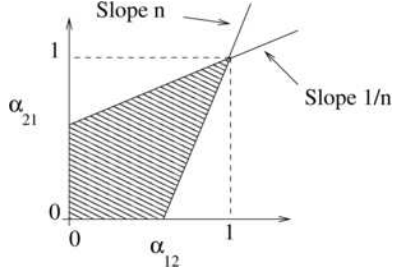

Fig. 1. *The persistence region of the Neuhauser–Pacala model.*

THEOREM 2.8 (Neuhauser and Pacala [13], Corollary 1). *Suppose that $\lambda = 1$. Write $n = |\mathcal{N}|$ for the number of lattice sites in a neighborhood. Species 1 competitively excludes species 2 if*

$$\alpha_{12} < \begin{cases} n\alpha_{21} - n + 1, & \text{for } \alpha_{21} \in \left(1 - \frac{1}{n}, 1\right], \\ \frac{1}{n}\alpha_{21} + 1 - \frac{1}{n}, & \text{for } \alpha_{21} > 1. \end{cases}$$

*Species 2 competitively excludes species 1 if*

$$\alpha_{12} > \begin{cases} \frac{1}{n}\alpha_{21} + 1 - \frac{1}{n}, & \text{for } \alpha_{21} \in (0, 1], \\ n\alpha_{21} - n + 1, & \text{for } \alpha_{21} > 1. \end{cases}$$

In particular, this result shows that the values of $(\alpha_{12}, \alpha_{21})$ for which both populations persist for all time are contained in the shaded region in Figure 1. This is a reduction from the range of values predicted by the mean field model. The case $\lambda = 1$ corresponds, in our setting, to taking $m_{ij} = m'_{ij}$, $\alpha = \alpha'$ and $M = M'$, so in view of Conjecture 2.3, we expect that the coexistence region for Model II in two dimensions corresponds to the whole region $[0, 1) \times [0, 1)$ in the $(\alpha_{12}, \alpha_{21})$-plane, that is, the region predicted by the mean field model.

Cox and Perkins [5] show that a sequence of processes following the Neuhauser–Pacala model, when suitably rescaled in space and time, converges to a superBrownian motion with a nontrivial drift. In low dimensions, they restrict to long-range models, whereas in dimensions $d \geq 3$, they can also consider the nearest-neighbor case. SuperBrownian motion has emerged as a universal limit of critical spatial systems above the critical dimension and these results can be seen as special cases of a general convergence theorem for perturbations of the voter model. In a recent preprint, Cox and Perkins [6] show that in dimensions $d \geq 3$, the drift in the superBrownian motion is connected to questions of coexistence in the Neuhauser–Pacala model. Using this connection, they obtain additional information about the parameter regions in which survival of one type (resp. coexistence) holds. The biologically relevant case $d = 2$ is a topic of their current research.

**3. $2N$-dependent oriented percolation.** We now turn to proving our results. Since our proofs will rely upon comparison with $2N$-dependent oriented percolation, we first briefly recall some well-known facts which can be found, for example, in Durrett [8]. The insistence on $2N$- instead of $N$-dependent percolation will be explained in the remark below Theorem 3.5.

Oriented percolation will be defined on the lattice

$$\mathcal{L} := \{(x,n) \in \mathbb{Z}^2 : x + n \text{ is even}, n \geq 0\}.$$

This set is made into a graph by inserting edges from $(x,n)$ to $(x+1, n+1)$ and to $(x-1, n+1)$. It is convenient to think of $n$ as time. We introduce a family of $\{0,1\}$-valued random variables $\omega(x,n)$ at sites $(x,n) \in \mathcal{L}$. A site $(x,n)$ is called *open* if $\omega(x,n) = 1$ and *closed* if $\omega(x,n) = 0$. Given such a family of random variables and integers $0 \leq m < n$, we say that $(y,n) \in \mathcal{L}$ can be reached from $(x,m)$ if there is a sequence of points $x = x_m, x_{m+1}, \ldots, x_n = y$ such that $|x_k - x_{k-1}| = 1$ and $\omega(x_k, k) = 1$ for $m \leq k \leq n$. We write this as $(x,m) \to (y,n)$. Finally, given an initial condition $\mathcal{W}_0 \subseteq 2\mathbb{Z} = \{x : (x,0) \in \mathcal{L}\}$, we may define a percolation process $\{\mathcal{W}_n\}_{n \geq 0}$ by setting, for each $n > 0$, $\mathcal{W}_n = \{y : (x,0) \to (y,n) \text{ for some } x \in \mathcal{W}_0\}$.

DEFINITION 3.1. Let $\theta \in (0,1)$ and $N \in \mathbb{N}$. We say that an oriented percolation process $\{\mathcal{W}_n\}_{n \geq 0}$, determined by $\{\omega(x,n)\}_{(x,n) \in \mathcal{L}}$, is $2N$-*dependent with density at least* $1 - \theta$ if, for any finite set of indices $I$ such that $\|(x_k, n_k) - (x_l, n_l)\| > 2N$ for all $k \neq l \in I$, we have $\mathbb{P}[\omega(x_k, n_k) = 0, k \in I] \leq \theta^{|I|}$.

Define $\mathcal{C}_0 = \{(y,n) \in \mathcal{L} : (0,0) \to (y,n)\}$ as the open cluster containing the origin. We say that *percolation occurs* if $|\mathcal{C}_0| = \infty$. We first cite a result which gives us a lower bound for the probability of percolation depending on $\theta$ and $N$. A proof can be found in Durrett [8].

THEOREM 3.2. *If* $\theta \leq 6^{-4(4N+1)^2}$, *then* $\mathbb{P}[|\mathcal{C}_0| < \infty] \leq 55\theta^{1/(4N+1)^2} \leq \frac{1}{20}$.

For particle system models of evolving populations, a standard strategy for showing survival is to construct a suitable coupling with oriented percolation. Our approach amounts to a modification of this strategy to cope with the interactions between the two populations, so we now describe the relevant comparison theorems. Once again, we are citing Durrett [8], but we also present a modified version of the results which are adequate for our purposes. In Durrett's terminology, let $\{\xi_i(n), i \in \mathbb{Z}^d\}_{n \geq 0}$ denote a translation invariant time-homogeneous finite-range flip process with state space $\Omega = \{0,1\}^{\mathbb{Z}^d}$, constructed from the usual graphical representation. Let $L \in \mathbb{N}$ be fixed. We write

$$(10) \qquad H = \{\xi \in \{0,1\}^{\mathbb{Z}^d} : \xi_i = 1 \ \forall i \in [-L/2, L/2]^d \cap \mathbb{Z}^d\}$$

and $m \cdot H$ for the translation of $H$ by some integer $m$ with respect to the first component, that is,

$$(11) \quad m \cdot H = \{\xi \in \{0,1\}^{\mathbb{Z}^d} : \xi_i = 1 \ \forall i \in mL\mathbf{e}_1 + [-L/2, L/2]^d \cap \mathbb{Z}^d\},$$

where $\mathbf{e}_1$ is the unit vector in the direction of the first component.

DEFINITION 3.3. Fix $N \in \mathbb{N}$ and $\theta \in (0,1)$. We shall say that the process $\{\xi(n)\}_{n \geq 0}$ fulfils the *classical comparison assumptions* for $N$ and $\theta$ if, for each configuration $\xi \in H$, there exists a "good event" $G_\xi$, measurable with respect to the graphical representation of the flip process inside $[-NL, NL]^d \times [0,1]$ and with $\mathbb{P}[G_\xi] > 1 - \theta$, so that if $\xi(0) = \xi$, then on $G_\xi$, we have $\xi(1) \in (+1) \cdot H \cap (-1) \cdot H$.

Unfortunately, the flip processes we are going to consider in the next section are functionals of more general underlying stochastic processes driven by independent Brownian motions and cannot be obtained from a graphical representation. However, a *modified* comparison result holds, based on the behavior of the Brownian increments. More explicitly, we are going to construct our flip process in terms of the system of stochastic differential equations $X, Y$ from Model I. Define $\{\bar{\xi}_i(n), i \in \mathbb{Z}^d\}_{n \in \mathbb{N}}$ by

$$(12) \qquad \bar{\xi}_i(n) = \begin{cases} 1, & \text{if } X_i(2n) > c_1 \text{ and } Y_j(2n) < c_2 \\ & \qquad \forall j \in i + [-bL, bL]^d \cap \mathbb{Z}^d, \\ 0, & \text{otherwise,} \end{cases}$$

where $c_1, c_2$ are finite positive numbers and $b$ and $L$ satisfy the assumption of Theorem 1.2. Note the different time scales for $\bar{\xi}$ and $X, Y$ (the usefulness of this time change will become clear in the next section) and observe that $X, Y$ are time-homogeneous and the underlying system of driving Brownian motions is translation invariant.

DEFINITION 3.4. Assume $N := (b+2)$ and $\theta \in (0,1)$. Define the events $H$ and $m \cdot H$ for some integer $m \in \mathbb{Z}$ in terms of the process $\{\bar{\xi}(n)\}_{n \geq 0}$ in the same way as in (10) and (11). Define

$$(13) \quad \mathcal{F}^*(NL, [0,2]) = \sigma\{B_i(s), B'_j(t) : s, t \in [0,2]; i, j \in [-NL, NL]^d \cap \mathbb{Z}^d\},$$

where $\{\{B_i(s)\}_{s \geq 0}, \{B'_j(t)\}_{t \geq 0}, i, j \in \mathbb{Z}^d\}$ is the family of independent standard Brownian motions as in Model I; see (2) and (3). We shall say that the process $\{\bar{\xi}(n)\}_{n \geq 0}$ fulfils the *modified comparison assumptions* for $NL$ and $\theta$ if, for each configuration $\bar{\xi} \in H$, there exists a "good event" $G_{\bar{\xi}} \in \mathcal{F}^*(NL, [0, 2])$, with $\mathbb{P}[G_{\bar{\xi}}] > 1 - \theta$, such that if $\bar{\xi}(0) = \bar{\xi}$, then on $G_{\bar{\xi}}$,

$$\bar{\xi}(1) \in (+1) \cdot H \cap (-1) \cdot H. \tag{14}$$

In other words, if $\bar{\xi}$ has all 1's in the box of side length $L$ about the origin at time 0, then at time 1 (measured in time units for the $\bar{\xi}$ process), with probability at least $1 - \theta$, it has successfully "invaded" the boxes of side length $L$ translated by $-L\mathbf{e}_1$ and $L\mathbf{e}_1$ in a way that is measurable with respect to the Brownian increments inside the box around the origin of side $2NL$ and up to time $[0, 2]$ (measured in time units for $X, Y$).

The following two classical theorems (see, e.g., Durrett [8]) apply in both of the above settings and complete this section.

THEOREM 3.5. *If the classical (resp. modified) comparison assumptions hold for $\xi$ (resp. $\bar{\xi}$) for some $N \in \mathbb{N}$ and $\theta \in (0, 1)$, we may define random variables $\omega(x, n)$ such that $\mathcal{X}_n := \{(m, n) \in \mathcal{L} : \xi(n) \in m \cdot H\}$ [resp. $\bar{\mathcal{X}}_n := \{(m, n) \in \mathcal{L} : \bar{\xi}(n) \in m \cdot H\}$] dominates a $2N$-dependent oriented percolation process $\{\mathcal{W}_n\}$, defined on $\mathcal{L}$, with initial configuration $\mathcal{W}_0 = \mathcal{X}_0$ (resp. $\bar{\mathcal{X}}_0$) and density parameter at least $1 - \theta$, that is, $\mathcal{W}_n \subseteq \mathcal{X}_n$ (resp. $\bar{\mathcal{X}}_n$) for all $n \in \mathbb{N}$, where $\mathcal{W}_n = \{y : (x, 0) \to (y, n) \text{ for some } x \in \mathcal{W}_0\}$.*

THEOREM 3.6. *Suppose that $\{\mathcal{W}_n\}_{n \geq 0}$ is a $2N$-dependent oriented percolation process, started from the trivial initial state $\mathcal{W}_0(x) = 1$ for all $x$. If $\theta \leq 6^{-4(4N+1)^2}$, then*

$$\liminf_{n \to \infty} \mathbb{P}[0 \in \mathcal{W}_{2n}] \geq \tfrac{19}{20}.$$

Eventually, this theorem will be the key to proving the coexistence (with positive probability) of $(X, Y)$ in Model I and Model II.

**4. Proofs.** The key to the proof of Theorem 1.2 is to consider the $Y$-population as providing a random environment in which the $X$-population evolves. Of course, the environment itself depends on the $X$-population, but we obtain some control of the behavior of the environment that is *independent* of the evolution of the $X$-population. This "decoupling" (and a symmetric argument for $Y$) then reduces the coexistence problem to that of survival of a single population: if the $X$- and $Y$-populations can each be shown to survive for all times with probability greater than one half, then longterm coexistence (with positive probability) will follow. We attack the question of survival of the $X$-population (resp. $Y$-population) by comparison with a $2N$-dependent oriented percolation process. To this end, we establish the existence of the corresponding "good events" as required in Definition 3.4.

4.1. *A spin system and estimation of related flip probabilities.* The main step is to construct two spin systems, one for each of the $X$- and $Y$-populations, that play the rôle of $\{\bar{\xi}_i(n), i \in \mathbb{Z}^d\}_{n\geq 0}$ of the last section for some suitable constants $c_1, c_2$. Indeed, we consider the spin system $\{\zeta_i(n), \eta_i(n), i \in \mathbb{Z}^d\}_{n\geq 0}$, where $\{\zeta_i(n), i \in \mathbb{Z}^d\}_{n\geq 0}$ is defined by

$$
(15) \qquad \zeta_i(n) = \begin{cases} 1, & \text{if } X_i(2n) > \dfrac{M}{K} \text{ and } Y_j(2n) < a'M' \\ & \qquad \forall j \in i + [-bL, bL]^d \cap \mathbb{Z}^d, \\ 0, & \text{otherwise}, \end{cases}
$$

and similarly,

$$
(16) \qquad \eta_i(n) = \begin{cases} 1, & \text{if } Y_i(2n) > \dfrac{M'}{K'} \text{ and } X_j(2n) < aM \\ & \qquad \forall j \in i + [-bL, bL]^d \cap \mathbb{Z}^d, \\ 0, & \text{otherwise}, \end{cases}
$$

where $K := 2\alpha M c + 1$, $K' := 2\alpha' M' c + 1$ and $a, a'$ are finite positive constants to be determined later [see (44) in the proof of Lemma 4.5 (resp. the symmetric result for $a$)]. Recall that $L$ denotes the range of the intraspecific interaction and $b$ denotes the smallest positive integer such that the range of $\{\gamma_{ij}\}$ (resp. $\{\gamma'_{ij}\}$) is less than $(b-1) \cdot L$.

With these definitions, one expects that if the system $\{\zeta(n), \eta(n)\}_{n\geq 0}$ exhibits longterm coexistence in discrete time, then the system $\{X(t), Y(t)\}_{t\geq 0}$ exhibits coexistence in continuous time and, in fact, this will follow from our proof. The convenience of the time change $n \mapsto 2n$ in (15) and (16) will become clear when carrying out the comparison arguments in Section 4.2.1.

*Outline of this subsection.* The notation from Section 4.1.1 is necessary to define the suitably measurable "good events" $G_\zeta, G_\eta$. Section 4.1.2 provides "flip probabilities" related to the spin system $\zeta$ via comparisons in terms of the aforementioned behavior of the one-dimensional diffusions, under the additional condition that the system $X$ evolves in a "safe environment," that is, given some bounds on the local $Y$-population. Finally, in Section 4.1.3, we will find conditions so that the "safe environment" assumption holds for the $Y$-population in a way that is independent of the evolution of the $X$-population, again by making use of comparisons to one-dimensional diffusions.

4.1.1. *Some notation and lattices of one-dimensional diffusions.*

DEFINITION 4.1. Let $\{\{B_i(s)\}_{s\geq 0}, \{B'_j(t)\}_{t\geq 0}, i, j \in \mathbb{Z}^d\}$ be the family of independent standard Brownian motions driving Model I. Fix $i \in \mathbb{Z}^d$. For $n \in \mathbb{N}, u > 0$, define

$$\mathcal{F}(i, n, u) := \sigma\{B_i(n+s) - B_i(n) : s \in [0, u]\}$$

and define $\mathcal{F}'(i, n, u)$ accordingly in terms of $B'$. Moreover, let

$$\mathcal{F}(i, NL, n, u) := \sigma\{B_j(n+s) - B_j(n) : s \in [0, u], j \in i + [-NL, NL]^d \cap \mathbb{Z}^d\},$$

and similarly define $\mathcal{F}'(i, NL, n, u)$ in terms of $B'$.

These $\sigma$-algebras will be used to construct suitably measurable events $G_{\bar{\zeta}}, G_{\bar{\eta}}$. Recall from Definition 3.4 that

$$\mathcal{F}^*(NL, [0, 2]) = \mathcal{F}(0, NL, 0, 2) \vee \mathcal{F}'(0, NL, 0, 2).$$

DEFINITION 4.2. Let $i \in \mathbb{Z}^d$ and assume that the constants $\alpha', \bar{M}, \lambda, \bar{U} > 0$ are chosen such that $\lambda \bar{U} > 2\bar{M}$. Moreover, let $D_1, D_2 > 0$. Then, for each $j \in i + [-NL, NL]^d \cap \mathbb{Z}^d$, define the one-dimensional diffusions $\{Z_j(t)\}_{t \geq 0}$, $\{\bar{Z}_j(t)\}_{t \geq 0}, \{\hat{Z}_j(t)\}_{t \geq 0}$ and $\{\tilde{Z}_j(t)\}_{t \geq 0}$, driven by independent standard Brownian motions $\{W_j(t)\}_{t \geq 0}$, by

(17) $$dZ_j(t) = \alpha'(\bar{M} - \lambda Z_j(t)) Z_j(t)\, dt + \sqrt{Z_j(t)}\, dW_j(t)$$

(logistic Feller diffusion),

(18) $$d\tilde{Z}_j(t) = D_1 \tilde{Z}_j(t)\, dt + \sqrt{\tilde{Z}_j(t)}\, dW_j(t)$$

(supercritical Feller diffusion),

(19) $$d\hat{Z}_j(t) = D_2\, dt + D_1 \hat{Z}_j(t)\, dt + \sqrt{\hat{Z}_j(t)}\, dW_j(t)$$

(supercritical Feller diffusion with constant positive immigration),

(20) $$d\bar{Z}_j(t) = \alpha'(\bar{M} - \lambda \bar{U}) \bar{Z}_j(t)\, dt + \sqrt{\bar{Z}_j(t)}\, dW_j(t)$$

(subcritical Feller diffusion).

Since each of the four diffusions admits a (continuous) unique strong solution, we may assume them to be driven by some given family of independent Brownian motions, in particular, those obtained either from (2) or from (3).

4.1.2. *Infection and recovery probabilities for the $X$-population.* Suppose that we are interested in the behavior of the $X$-population within the time interval $[n, n+1]$ at site $i$ and that we already know [recall (15)] that

(21) $$\max_{t \in [n, n+1]} Y_j(t) < 2a'M' \qquad \forall j \in i + [-bL, bL]^d \cap \mathbb{Z}^d.$$

Assume the interspecific competition $\{\gamma_{ij}\}$ is chosen such that (21) implies

(22) $$\max_{t \in [n, n+1]} \sum_{l \in \mathbb{Z}^d} \gamma_{jl} Y_l(t) < 1 \qquad \forall j \in i + [-L/2, L/2]^d \cap \mathbb{Z}^d.$$

This is possible since the range of $\{\gamma_{jl}\}$ is, by assumption, less than $(b-1)L$ [choose, e.g., $\sum_j \gamma_{ij} < (2a'M')^{-1}$]. We will later (in Section 4.1.3) construct events that are measurable with respect to either $\mathcal{F}'(i, NL, n, 2)$ or $\mathcal{F}'(i, NL, n-1, 2)$, which imply (21) and are of sufficiently high probability [cf. (39) and (40)]. For the moment, to aid intuition and to simplify notation, we will say that, in either case, a suitably measurable "*safe environment condition* $G'_{\mathrm{sec}}(i, n)$ holds at site $i$ and time $n$," which implies that (21) [and for sufficiently small $\gamma_{ij}$, also (22)] holds and which will be explicitly determined later.

We now consider "flip probabilities" for the $X$-population that are closely linked to the flip probabilities of the $\zeta$-population, introduced in (15) and (16), under the "safe environment condition" $G'_{\mathrm{sec}}(i, n)$ at site $i$ and at time $n \in \mathbb{N}$.

LEMMA 4.3 (Infection and nonrecovery). *Let $n \in \mathbb{N}$ and $i \in \mathbb{Z}^d$. Let $\alpha, \{m_{ij}\}, \{\lambda_{ij}\}$ be fixed. Given the parameters for the $Y$-population and some $a' > 0$, choose $\{\gamma_{ij}\}$ such that* (21) *implies* (22). *Then, for any $\varepsilon \in (0,1)$, there exists a finite constant $M_0 > 0$ such that if $M > M_0$ and $K := 2\alpha Mc + 1$, for each $j \in i + [-L/2, L/2]^d \cap \mathbb{Z}^d$, there exist events*

(23) $\quad G_{\mathrm{nonrec}}(i, n) \in \mathcal{F}(i, n, 1) \quad$ and $\quad G_{\mathrm{infec}}(i, j, n) \in \mathcal{F}(i, L/2, n, 1),$

*both measurable w.r.t. $\mathcal{F}(i, NL, n, 1)$, such that the following holds*:
(i) *we have*

$$\{\{X_i(n) > M/K\} \cap G'_{\mathrm{sec}}(i, n) \cap G_{\mathrm{nonrec}}(i, n)\} \subset \{X_i(n+1) > M/K\}$$

*("nonrecovery") and for the "nonrecovery probability" $p_{\mathrm{nonrec}}(i, n)$, we have the bound*

(24) $$p_{\mathrm{nonrec}}(i, n) := \mathbb{P}[G_{\mathrm{nonrec}}(i, n)] > 1 - \varepsilon;$$

(ii) *moreover*,

$$\{\{X_i(n) \leq M/K\} \cap \{\exists j : m_{ij} > 0, X_j(n) > M/K\} \cap G'_{\mathrm{sec}}(i, n)\}$$
$$\subset \{X_i(n+1) > M/K\}$$

*("infection by an occupied neighbor") and for the "infection probability" $p_{\mathrm{infec}}(i, j, n)$, we have the bound*

(25) $$p_{\mathrm{infec}}(i, j, n) := \mathbb{P}[G_{\mathrm{infec}}(i, j, n)] > 1 - \varepsilon.$$

PROOF. (i) We distinguish the two cases $X_i(n) \in (M/K, (3/2)M/K)$ and $X_i(n) \geq (3/2)M/K$.

*Case* 1. Suppose that $X_i(n) \geq (3/2)M/K$ and introduce the first hitting time of level $M/K$ from above after time $n$:

(26) $$\tau^{X_i}_{M/K}(n) := \inf\{t > n : X_i(t) = M/K\}.$$

Our goal is to establish the existence of a suitably measurable event $G_{\text{nonrec}}(i,n) \in \mathcal{F}(i,n,1)$ so that $G_{\text{nonrec}}(i,n)$ implies, under the above conditions, that $\tau^{X_i}_{M/K}(n) > 1$. To this end, we set up a suitable comparison to a one-dimensional diffusion. Indeed, rearranging the drift in equation (2), as long as $Y_j(t) < a'M'$ for all $j \in i + [-bL, bL]^d \cap \mathbb{Z}^d$ and hence

$$\sum_l \gamma_{jl} Y_l(t) < 1 \qquad \text{for all } j \in i + [-L/2, L/2]^d \cap \mathbb{Z}^d$$

holds, and as long as $X_i \leq 2M/K$, we have

(27) $$dX_i(t) \geq \sum_{j \in \mathbb{Z}^d} \left( m_{ij} - \alpha \frac{2M}{K} \lambda_{ij} \right) X_j(t)\, dt$$
$$+ \left( \alpha M - \sum_{j \in \mathbb{Z}^d} m_{ij} - \alpha \right) X_i(t)\, dt + \sqrt{X_i(t)}\, dB_i(t).$$

We now check that the first component of the drift on the right-hand side is positive. Indeed, from the assumption (4), we obtain

(28) $$m_{ij} - \alpha \frac{2M}{K} \lambda_{ij} > m_{ij} - \alpha \frac{2M}{K} c m_{ij},$$

which is positive by our choice of $K = 2\alpha M c + 1$, for all $j \in \mathbb{Z}^d$. Moreover, we have, for each $M > 1$, that $M/K = M/(2\alpha M c + 1) \in ((2\alpha c + 1)^{-1}, (2\alpha c)^{-1})$. Under these conditions, (27) implies that

(29) $$dX_i(t) \geq \left( \alpha M - \sum_{j \in \mathbb{Z}^d} m_{ij} - \alpha \right) X_i(t)\, dt + \sqrt{X_i(t)}\, dB_i(t)$$

and so while $X_i \in [0, 2M/K]$, using Corollary A.2 to the Ikeda–Watanabe Comparison Theorem A.1 (both found in the Appendix), we may compare $X_i$ to a dominated supercritical Feller diffusion $\tilde{Z}_i$ defined in (18), with initial value $\tilde{Z}_i(n) := (3/2)M/K$ and $D_1 = D_1(M) = (\alpha M - \sum_{j \in \mathbb{Z}^d} m_{ij} - \alpha)$, and driven by the same Brownian motion, that is, $\{W_i(t)\}_{t\geq 0} := \{B_i(t)\}_{t\geq 0}$. It is an important observation in Corollary A.2 that actually more is true: our domination argument holds not only up to the time when $X_i$ leaves the interval $[0, 2M/K]$ for the first time, but, in fact, for as long as $\tilde{Z}_i$ takes values inside this interval, that is, up to the first exit time $\tau^{\tilde{Z}_i}_{2M/K}$ [defined as in (26)].

Note that for $M > \sum_j m_{ij}/\alpha$, the "supercriticality" (i.e., positive drift) $D_1 = (\alpha M - \sum m_{ij} - \alpha)$ in (29) tends to $\infty$ as $M \to \infty$, while maintaining the condition $X_i(n) \geq 1/(2\alpha c)$.

We now make use of the comparison. Indeed, for $t \geq n$, as long as $\tilde{Z}_i(t)$ stays inside the interval $[0, 2M/K]$ and given that initially $X_i(n) \geq \tilde{Z}_i(n) := (3/2)M/K$, we have that $X_i$ dominates $\tilde{Z}_i$. To obtain a comparison that is valid throughout the whole time interval $[n, n+1]$, we go one step further and modify $\tilde{Z}_i$ so that whenever $\tilde{Z}_i$ hits level $2M/K$ (and thus is about to leave the area in which the comparison holds true), we restart the process $\tilde{Z}_i$ at level $(3/2)M/K$ and repeat this procedure as often as necessary, so that the comparison holds for all times $t \in [n, n+1]$. More precisely, we define a sequence of stopping times, beginning with $\nu^{\tilde{Z}_i}_{2M/K}(n, 1) := \tau^{\tilde{Z}_i}_{2M/K}(n)$, restart the $\tilde{Z}_i$ process at this time, setting $\tilde{Z}_i(\nu^{\tilde{Z}_i}_{2M/K}(n, 1)) := (3/2)M/K$, and then iterate this procedure, considering, for $m \in \mathbb{N}$,

$$(30) \quad \nu^{\tilde{Z}_i}_{2M/K}(n, m+1) := \inf\{t > \nu^{\tilde{Z}_i}_{2M/K}(n, m) : \tilde{Z}_i(t) = (2M)/K\},$$

and again restarting the $\tilde{Z}_i$ process accordingly, that is, setting

$$\tilde{Z}_i(\nu^{\tilde{Z}_i}_{2M/K}(n, m+1)) = (3/2)(M/K).$$

Note that $\nu^{\tilde{Z}_i}_{2M/K}(n, m) \uparrow \infty$ a.s. as $m \to \infty$. For definiteness, set $\nu^{\tilde{Z}_i}_{2M/K}(n, 0) := n$. We define the i.i.d. positive lengths of the corresponding upcrossing intervals from $(3/2)M/K$ to $2M/K$ for $m \geq 1$ by

$$(31) \quad \tilde{T}_m := \nu^{\tilde{Z}_i}_{2M/K}(n, m) - \nu^{\tilde{Z}_i}_{2M/K}(n, m-1).$$

Now observe that there is an event $G^1_{\text{nonrec}}(i, n)$, defined only in terms of the Brownian increments $\{B_i(n+s) - B_i(n) : s \in [0, 1]\}$, and hence being an element of $\mathcal{F}(i, n, 1)$, such that if we start our modified diffusion $\tilde{Z}_i$ in $\tilde{Z}_i(n) = (3/2)M/K$, this event $G^1_{\text{nonrec}}(i, n)$ actually equals $\{\tau^{\tilde{Z}_i}_{M/K}(n) > n+1\}$. [The set $G^1_{\text{nonrec}}(i, n)$ contains all such $\omega$, such that the corresponding Brownian increments lead to the desired behavior if they drive the modified diffusion $\tilde{Z}_i$ started at time $n$ in $(3/2)M/K$.] Moreover, by our comparison, the event $\{\tau^{\tilde{Z}_i}_{M/K}(n) > n+1\}$ implies that $\{\tau^X_{M/K}(n) > n+1\}$, which, in turn, implies that $X_i(n+1) > M/K$.

It remains to show that the event $G^1_{\text{nonrec}}(i, n)$ has sufficiently high probability. To this end, note that the number of upcrossings of the modified and suitably restarted process $\tilde{Z}_i$ from level $(3/2)M/K$ to level $2M/K$ before the first downcrossing from $(3/2)M/K$ to $M/K$ is a geometric random variable with positive parameter $\tilde{q}_{M,K} := \mathbb{P}[\tau^{\tilde{Z}}_{M/K}(n) < \tau^{\tilde{Z}}_{2M/K}(n)]$. By standard speed measure and scale function computations, it is possible to show, for sufficiently large $D_1$ and hence $M$, that $\tilde{q}_{M,K} \leq \exp(-D_1 M/K)$; see Blath, Meredith and Etheridge [1] for details. Moreover, by a similar computation, the expected time for $\tilde{Z}_i$ to exit from the interval $(M/K, 2M/K)$ when being started in $(3/2)M/K$ is bounded below by $\frac{1}{32D_1}$ for sufficiently large $D_1$. Hence, for such $D_1$,

$$\mathbb{E}[T_1] \geq \mathbb{E}[\tau^{\tilde{Z}_i}_{(M/K, 2M/K)}] \geq \frac{1}{32D_1}. \tag{32}$$

Now let $\tilde{\mathcal{D}}$ denote the number of upcrossings before the first "success," that is, a downcrossing from $(3/2)M/K$ to $M/K$. For each $\tilde{N} \in \mathbb{N}$, we may then write

$$\mathbb{P}[\tau^{\tilde{Z}_i}_{M/K}(n) \leq n+1] = \mathbb{P}[\tau^{\tilde{Z}_i}_{M/K}(n) \leq n+1; \tilde{\mathcal{D}} < \tilde{N}]$$
$$+ \mathbb{P}[\tau^{\tilde{Z}_i}_{M/K}(n) \leq n+1; \tilde{\mathcal{D}} \geq \tilde{N}]$$
$$\leq \mathbb{P}[\tilde{\mathcal{D}} < \tilde{N}] + \mathbb{P}\left[\sum_{i=1}^{\tilde{\mathcal{D}}} \tilde{T}_i < 1; \tilde{\mathcal{D}} \geq \tilde{N}\right]$$
$$\leq \tilde{N} \exp(-D_1 M/K) + \mathbb{P}\left[\sum_{i=1}^{\tilde{N}} T_i < 1\right],$$

using Bernoulli's inequality. Since, by (32), for large $D_1$, the expectation of the length of the i.i.d. upcrossing intervals $\{\tilde{T}_i\}$ of the modified and suitably restarted process is bounded below by $1/(32D_1)$, the number of such upcrossing intervals up to time 1 is at most of order $D_1$. Hence, by the Law of Large Numbers, we can find a constant $\tilde{d}$ such that for $\tilde{N} := \tilde{d} \cdot D_1$ and all sufficiently large $D_1$, the last term on the right-hand side is bounded by $\varepsilon/4$. Since the first term on the right-hand side still decreases exponentially in $D_1$ once $\tilde{d}$ is fixed (the linearly increasing prefactor being squashed), for $D_1$ and hence $M$ sufficiently large, this bound holds simultaneously for the first and the last term and we arrive at the desired result: under the above conditions, with $\tilde{Z}_i(n) = (3/2)M/K$, we have $\mathbb{P}[\tau^{\tilde{Z}_i}_{M/K}(n) \leq n+1] \leq \varepsilon/2$ for $D_1$ and hence $M$ sufficiently large, which in turn implies $\mathbb{P}[G^1_{\text{nonrec}}(i,n)] > 1 - \frac{\varepsilon}{2}$, so that Case 1 of part (i) follows.

*Case* 2. Now suppose that $M/K < X_i(n) < (3/2)M/K$. In this case, we cannot find a uniform lower bound on the probability on the previously considered event $\{\tau^{\tilde{Z}_i}_{M/K}(n) > n+1\}$, and hence on the probability of $\{\tau^{X_i}_{M/K}(n) > n+1\}$, that is sufficiently large. However, we may still use the same comparison as above to a dominated supercritical Feller diffusion $\tilde{Z}_i$, so that the comparison works as long as $\tilde{Z}_i$ stays below $2M/K$. This time, set $\tilde{Z}_i(n) = M/K < X_i(n)$ and observe that there is a constant $M_0^2 > 0$ such that for all $M > M_0^2$, the deterministic drift in the supercritical Feller diffusion $\tilde{Z}_i$ will achieve two goals with sufficiently high probability: first, make $\tilde{Z}_i$ hit level $(3/2)M/K$ within the time interval $[n, n+1/2]$ with sufficiently high probability and second, after hitting level $(3/2)M/K$, arguing just as in the first part of the lemma, ensure that there will be no further downcrossing from $(3/2)M/K$ to $M/K$ up to time $[n+1]$. Thus, once again, we can find a measurable event $G_{\text{nonrec}}^2(i,n) \in \mathcal{F}(i,n,1)$, depending only on the corresponding Brownian increments, so that given $\tilde{Z}_i(n) = M/K$, by comparison, $G_{\text{nonrec}}^2(i,n)$ implies that $X_i(n+1) > M/K$ and, moreover, that $\mathbb{P}[G_{\text{nonrec}}^2(i,n)] > 1 - \varepsilon/2$. Hence, the result also holds in Case 2. Finally, in view of both cases, choose

$$G_{\text{nonrec}}(i,n) := G_{\text{nonrec}}^1(i,n) \cap G_{\text{nonrec}}^2(i,n) \in \mathcal{F}(i,n,1) \subset \mathcal{F}(i,NL,n,1)$$

and part (i) follows.

To prove part (ii), we begin with some preliminary considerations. Note that by using the same comparison and similar arguments as before, again considering suitable up- and downcrossings, this time from $M/K$ down to $M/(2K)$, we can actually go one step further and find a finite constant $M_0^3$ such that if $M > M_0^3$, and for $j \in \mathbb{Z}^d$ such that $m_{ij} > 0$, there exists an event $G_{\text{per-occ}}(i,j,n) \in \mathcal{F}(i,NL,n,1)$, such that given $X_j(n) > M/K$ and $G'_{\text{sec}}(i,n)$, the event $G_{\text{per-occ}}(i,j,n)$ implies $X_j(n+1) > M/K$ and $\tau_{M/(2K)}^{X_j}(n) > n+1$ and, moreover, we have

$$\mathbb{P}[G_{\text{per-occ}}(i,j,n)] > 1 - \varepsilon/2. \tag{33}$$

Note that, once again, we use the assumption that the range of the $\{\gamma_{ij}\}$ is less than $(b-1)L$ so that the "safe environment condition," in particular (22), allows comparisons of the above type also for site $j$.

We are now prepared to consider the *infection probability* at a site $i$ in the presence of at least one occupied neighbor, say, at $j^*$. Again, assuming (4), we use a comparison based on Corollary A.2. This theorem can be found in the Appendix. This time, we rearrange the drift so as to highlight the rôle of immigration of mass to an unoccupied site from occupied neighbors. Once immigrated, we can then compare the evolution of the mass to a supercritical continuous-state branching process, as before. Indeed, considering the drift in equation (2), observe that as long as $X_i(t) \le 2M/K$ and given the existence of at least one neighbor at some site $j^* \in \mathbb{Z}^d$ with $m_{ij^*} > 0$ and $X_{j^*}(n) > M/K$ (noting that $m_{ij^*}$ is bounded below by some $\delta > 0$ since the family $\{m_{ij}\}$ is of finite range), we have that, as long as $t$ satisfies

$$n \le t < \tau_{M/(2K)}^{X_{j^*}}(n), \tag{34}$$

by our choice $K = 2\alpha Mc + 1$,

$$
\begin{aligned}
dX_i(t) &\geq \frac{m_{ij^*} - (1/c)\lambda_{ij^*}}{4\alpha c + 2}\, dt \\
&\quad + \left(\alpha M - \sum_{j \in \mathbb{Z}^d} m_{ij} - \alpha\right) X_i(t)\, dt + \sqrt{X_i(t)}\, dB_i(t),
\end{aligned}
\tag{35}
$$

assuming $M > \max\{1, \sum_j m_{ij}/\alpha\}$. Thus, at the uninfected site $i$, after time $n$, as long as (34) holds, we may compare the evolution of the process $X_i$ to a dominated supercritical branching process $\hat{Z}_i$ with constant immigration, as defined in (19), and driven by the same Brownian motion, that is, $\{W_i(t)\}_{t \geq 0} := \{B_i(t)\}_{t \geq 0}$, where $D_2 = (m_{ij^*} - \lambda_{ij^*}/c)(4\alpha c + 2)^{-1} > 0$ and $D_1 = (\alpha M - \sum_j m_{ij} - \alpha) > 0$. Note that we here use the fact that we have assumed the strict inequality $\frac{1}{c}\lambda_{ij^*} < m_{ij^*}$ from (4) to obtain a strictly positive immigration rate. Also, note that the rate of immigration is bounded below independently of $M$. Again, the "supercriticality" $(\alpha M - \sum_j m_{ij} - \alpha)$ tends to $\infty$ as $M$ tends to $\infty$. Hence, arguing as before in part (i), this time starting the dominated process in $\hat{Z}_i(n) = 0 \leq X_i(n)$, stopping $\hat{Z}_i$ once it reaches level $2M/K$ and then restarting at $(3/2)M/K$ if necessary, we may find a constant $M_0^4 > 0$ such that if $M > M_0^4$, under the above conditions, the event

$$\hat{Z}_i(n+1) > M/K \qquad \text{implies } X_i(n+1) > M/K$$

and has probability greater than $\varepsilon/2$. As before, there is an event $G^*_{\text{infec}}(i, j^*, n) \in \mathcal{F}(i, n, 1)$, depending solely on the Brownian increments at $i$ within the time interval $[n, n+1]$, so that if we start our modified diffusion at time $n$ in $\hat{Z}_i(n) = 0$, driven by the corresponding Brownian increments, we have $G^*_{\text{infec}}(i, j^*, n) = \{\hat{Z}_i(n+1) > M/K\}$. Moreover,

$$\mathbb{P}[G^*_{\text{infec}}(i, j^*, n)] > 1 - \varepsilon/2.$$

Now, observe that, due to our preparations, (34) is actually being guaranteed by $G_{\text{per-occ}}(i, j^*, n)$ up to time $n+1$. Combining both events (which are actually independent), we define

$$G_{\text{infec}}(i, j^*, n) := G^*_{\text{infec}}(i, j^*, n) \cap G_{\text{per-occ}}(i, j^*, n) \in \mathcal{F}(i, L/2, n, 1)$$

and finally see that $\mathbb{P}[G_{\text{infec}}(i, j^*, n)] > 1 - \varepsilon$. Together with the fact that given $X_{j^*}(n) > M/K$, our event $G_{\text{infec}}(i, j, n)$ implies $X_i(n+1) > M/K$, defining $M_0 = \max\{M_0^1, M_0^2, M_0^3, M_0^4\}$ completes the proof of Lemma 4.3. □

4.1.3. *Control of the environment.* We now find estimates for suitably measurable events, a combination of which will later provide the "safe environment condition" $G'_{\text{sec}}(i,n)$ at some site $i \in \mathbb{Z}^d$ and time $n \in N$, which will, in turn, imply (21). Again, this is done via suitable comparisons to one-dimensional diffusions.

LEMMA 4.4 (Control of the environment). *Let $n \in \mathbb{N}$ and $i \in \mathbb{Z}^d$. Let $\alpha', \{m'_{ij}\}, \{\lambda_{ij}\}$ be fixed. Then, for any $\varepsilon \in (0,1)$, there is a finite constant $v'_0 > 0$ such that for all $v' > v'_0$, there exists an event $E'(i,n)[v'] \in \mathcal{F}'(i,n,2)$ such that for all $\alpha' M' > \sum_j m'_{ij}$,*

$$\mathbb{P}[E'(i,n)[v']] > 1 - \varepsilon.$$

*Moreover,*

$$(36) \qquad E'(i,n)[v'] \subset \left\{ \sup_{0 \leq s \leq 1} Y_i(n+1+s) < 2v'M' \right\}$$

*and*

$$(37) \qquad \{E'(i,n)[v'] \cap \{Y_i(n) \in [0, v'M']\}\} \subset \left\{ \sup_{0 \leq s \leq 2} Y_i(n+s) < 2v'M' \right\}.$$

*These results hold true for any choice of $\{\gamma'_{ij}\}$ and any parameter values for the $X$-population in Model* I.

REMARKS. 1. In the proof, we use the assumption that there is a constant $c > 0$ such that $m'_{ij} < c\lambda'_{ij}$. We do not believe that this is necessary.

2. That the result should be true is again due to the fact that the downward drift resulting from overcrowding in a site is quadratic, whereas the upward drift due to reproduction in the population is only linear. Moreover, for sufficiently crowded sites, immigration from neighboring sites is being compensated by intraspecific competition.

SKETCH OF PROOF OF LEMMA 4.4. Again, the proof relies on a suitable comparison to one-dimensional diffusions. This time, the comparison is set up in a way such that $Y_i$ is being *dominated* by a one-dimensional diffusion. Indeed, notice that, for $t \geq n$, as long as $Y_i(t) > m'_{ij}/(\alpha' \lambda'_{ij})$ for all $j \neq i$, an informal calculation shows that

$$(38) \qquad dY_i(t) \leq \alpha' \left( M' - \frac{1}{\alpha'} \sum_{j \neq i, j \in \mathbb{Z}^d} m'_{ij} - \lambda'_{ii} Y_i(t) \right) Y_i(t)\, dt + \sqrt{Y_i(t)}\, dB'_i(t).$$

Hence, the immigration of mass from site $j$ is compensated for by the downward drift due to crowding at site $j$ and we may compare the evolution of $Y_i$ (again applying the Corollary A.2 to the Ikeda–Watanabe Comparison theorem) to that of the solution $Z_i$ of the logistic Feller diffusion (17), this time with driving Brownian motion given by $\{W_i(t)\}_{t\geq 0} := \{B'_i(t)\}_{t\geq 0}$. It is then possible, in a way similar to the proof of Lemma 4.3, to construct the suitably measurable events $E'(i,n)[v']$ in terms of $Z_i$ with the required properties. Here, the quadratic downward drift in the logistic Feller diffusion does the trick, as standard speed-measure and scale-function computations for $Z_i$ show. Full details of this proof can again be found in Blath, Etheridge and Meredith [1]. □

REMARK. Note that for $a' > v'_0$ in Lemma 4.4 (and, of course, also for the stronger condition $a'/2 > v'_0$), at time $n \in \mathbb{N}$, recalling $N = b + 2$,

$$(39) \quad \bigcap_{j \in i+[-bL,bL]^d \cap \mathbb{Z}^d} E'(j,n-1)[a'] \in \mathcal{F}'(i,bL,n-1,2) \subset \mathcal{F}'(i,NL,n-1,2)$$

implies (21). Moreover, together with the additional condition $Y_j(n) \leq a'M'$ for all $j \in i + [-bL, bL]^d \cap \mathbb{Z}^d$, both

$$(40) \quad \bigcap_{j \in i+[-bL,bL]^d \cap \mathbb{Z}^d} E'(j,n-1)[a'] \in \mathcal{F}'(i,bL,n-1,2) \subset \mathcal{F}'(i,NL,n-1,2)$$

and

$$(41) \quad \bigcap_{j \in i+[-bL,bL]^d \cap \mathbb{Z}^d} E'(j,n)[a'] \in \mathcal{F}'(i,bL,n,2) \subset \mathcal{F}'(i,NL,n,2)$$

imply (21), too. Thus, all three events can be considered as instances of the "safe environment condition" $G'_{\text{sec}}(i,n)$ in the sense of Lemma 4.3.

4.2. *Comparison arguments.*

4.2.1. *Comparison of $\zeta$ (resp. $\eta$) to $2N = 2(b+2)$-dependent oriented percolation.* We first focus on the $X$- (resp. $\zeta$-) population and construct a suitable "good event" $G_\zeta$. Recall that by our technical assumption, $L$ denotes the maximum range of both migration matrices $\{m_{ij}\}, \{m'_{ij}\}$ and $b$ is the smallest positive integer such that the range of $\{\gamma_{ij}\}$ and $\{\gamma'_{ij}\}$ is less than $(b-1) \cdot L$. Also, recall from (15) the definition

$$(42) \quad \zeta_i(n) = \begin{cases} 1, & \text{if } X_i(2n) > \dfrac{M}{K} \text{ and } Y_j(2n) < a'M' \\ & \forall j \in i + [-bL, bL]^d \cap \mathbb{Z}^d, \\ 0, & \text{otherwise} \end{cases}$$

and from the comparison assumptions (10) the notation $\zeta(n) \in H$ if $\zeta_i(n) = 1$ for all $i \in [-L/2, L/2]^d \cap \mathbb{Z}^d$. Finally, recall from (11) the notion of translation of $H$ by $mL$, for some $m \in \mathbb{Z}$, denoted by $m \cdot H$.

LEMMA 4.5 (Comparison). *Let $\theta \in (0,1)$. Suppose we are given fixed parameters $\alpha, \alpha', \{m_{ij}\}, \{m'_{ij}\}, \{\lambda_{ij}\}, \{\lambda'_{ij}\}$. Then, under the assumptions on Model I for Theorem 1.2, there exist finite constants $M_0 > 0$ and $a' > 0$ such that if $M > M_0$, $K := 2\alpha Mc + 1$ and $M' > \sum_j m'_{ij}/\alpha'$, then there is a finite $\gamma = \gamma(a'M') > 0$ such that if $\sum_j \gamma_{ij} < \gamma$, and for all $\{\gamma'_{ij}\}$, the process $\{\zeta(n)\}_{n \geq 0}$ fulfils the* modified comparison assumptions (14) *for $NL$ and $\theta$. In particular, for each configuration $\zeta \in H$, there exists a "good event"*

$$G_\zeta \in \mathcal{F}^*(NL, [0,2]),$$

*where $N = b + 2$, with $\mathbb{P}[G_\zeta] > 1 - \theta$, such that if $\zeta(0) = \zeta$, then on $G_\zeta$,*

$$\zeta(1) \in (+1)H \cap (-1)H.$$

*Consequently, using Theorem 3.5, the process $\mathcal{X}_n := \{(m,n) \in \mathcal{L} : \zeta(n) \in m \cdot H\}$ dominates a $2N$-dependent oriented percolation process $\{\mathcal{W}_n\}_{n \geq 0}$ on $\mathcal{L}$ with density at least $1 - \theta$ and initial condition $\mathcal{W}_0 = \mathcal{X}_0$.*

REMARKS. 1. A similar result is true for the $Y$- (resp. $\eta$-) population which, given $\theta > 0$, produces a similar threshold $M'_0$ and parameters $M', a, \gamma'$, which allow a comparison to a $2N = 2(b+2)$-dependent oriented percolation process of density at least $1 - \theta$ via a similar "good event" $G_\eta$.

2. The available degree of freedom in the choice of $\{\gamma'_{ij}\}$ in this result is crucial for the simultaneous comparison of $\{\zeta(n)\}_{n \geq 0}$ and $\{\eta(n)\}_{n \geq 0}$ which we will need to consider later. It is due to the fact that our results for the "control of the environment" in Lemma 4.4 are entirely independent of these $\{\gamma'_{ij}\}$, since competition by $X$ only facilitates the "good environment condition" determined in terms of $Y$.

PROOF OF LEMMA 4.5. Fix $\theta > 0$ and let

$$(43) \qquad \varepsilon = \frac{1}{2} \frac{\theta}{(4(b+2)L)^d}.$$

We begin with the specification of consistent parameter values for our model that will lead to the required comparison. First, note that all of the constants $\alpha, \alpha', \{m_{ij}\}, \{m'_{ij}\}, \{\lambda_{ij}\}, \{\lambda'_{ij}\}, b, N, L$ will remain fixed throughout what follows. The only values we need to adjust suitably in order to produce the proof are $M, M', a', \{\gamma_{ij}\}$. The proof is entirely independent of the choice of $\{\gamma'_{ij}\}$ (provided all parameter values remain compatible with the assumptions of Theorem 1.2).

First, we choose sufficiently large $a'$ so that for any $M' > \sum_j m'_{ij}/\alpha'$,

$$(44) \qquad a' > \frac{1}{2\alpha' M'c + 1}$$

and, moreover, such that $a'/2 > v'_0$ in Lemma 4.4 with the above $\varepsilon$. Then, for each $i$ and $n$, we have the bound

(45) $$\mathbb{P}[E'(i,n)[a'/2]] > 1 - \varepsilon.$$

Note that this bound does not depend on $\{\gamma'_{ij}\}$ (and obviously not on $\{\gamma_{ij}\}$). From now on, $a'$ remains fixed. Define, for any $M' > \sum_j m'_{ij}/\alpha'$, the constant $\gamma = \gamma(a', M') := (2a'M')^{-1}$ so that, for each $i$, we have $\sum_{j \in \mathbb{Z}^d} \gamma_{ij} 2a'M' < 1$.

Finally, we can find $M_0 > 0$ such that for all $M > M_0$, the bounds of Lemma 4.3 for the "infection" and "nonrecovery probabilities" hold with our choice of $\varepsilon$.

We now check that with these parameter values for $M_0, a', \gamma$ and for all $M > M_0$, assuming $\zeta(0) \in H$, there is a "good event" $G_\zeta \in \mathcal{F}^*(NL, [0, 2])$, which implies

(46) $$\zeta(1) \in (+1)H \cap (-1)H$$

and has probability at least $1 - \theta$. Recall that $\zeta(0) \in H$ means:

- $X_i(0) > M/K$, where $K = 2\alpha Mc + 1$, for all $i \in [-L/2, L/2]^d \cap \mathbb{Z}^d$,
- $Y_j(0) < a'M'$, for all $j \in i + [-bL, bL]^d \cap \mathbb{Z}^d, i \in [-L/2, L/2]^d \cap \mathbb{Z}^d$.

To construct $G_\zeta$, recall that one time step for $\zeta$ corresponds to two time units for $X$ and $Y$. We split the corresponding time interval $[0, 2]$ into two parts, $[0, 1]$ and $[1, 2]$. By Lemma 4.4, applied with $v' = a' > v'_0$, we see that

$$\mathbb{P}\left[\bigcap_i E'(i, 0)[a'] : i \in [-(b+1/2)L, (b+1/2)L]^d \cap \mathbb{Z}^d\right]$$
$$> 1 - (2(b+1/2)L)^d \varepsilon,$$

and recall that this event, denoted by $E'(0, (b+1/2)L, 0, 2)[a'] \in \mathcal{F}'(0, (b+1/2)L, 0, 2)$ for short, implies, since $\zeta(0) \in H$, by Lemma 4.4, that

$$\sup_{0 \le s \le 2} \{Y_i(s) : i \in [-(b+1/2)L, (b+1/2)L]^d \cap \mathbb{Z}^d\} < 2a'M'.$$

Next, Lemma 4.3 tells us that [recall $G_{\text{nonrec}}(i, 0) \in \mathcal{F}(i, 0, 1)$ for $n = 0$],

$$\mathbb{P}\left[\bigcap_i G_{\text{nonrec}}(i, 0) : i \in [-L/2, L/2]^d \cap \mathbb{Z}^d\right] > 1 - L^d \varepsilon$$

and this event, denoted by $G_{\text{nonrec}}(0, L/2, 0, 1) \in \mathcal{F}(0, L/2, 0, 1)$, implies

$$X_i(1) > \frac{M}{2\alpha Mc + 1} \quad \text{for all } i \in [-L/2, L/2]^d \cap \mathbb{Z}^d.$$

Now, applying Lemma 4.4 once again, this time with $v' = a'/2$, we see that

$$\mathbb{P}\left[\bigcap_i E'(i, 0)[a'/2] : i \in [-NL, NL]^d \cap \mathbb{Z}^d\right] > 1 - \varepsilon(2NL)^d$$

and this event, denoted by $E'(0, NL, 0, 2)[a'/2] \in \mathcal{F}'(0, NL, 0, 2)$, ensures that
$$\sup_{s \in [1,2]} \{Y_i(s) : i \in [-NL, NL]^d \cap \mathbb{Z}^d\} < a'M'.$$
Combining all of these events, we have guaranteed that with probability at least $1 - (2NL)^d \varepsilon - L^d \varepsilon - (2(b+1/2)L)^d \varepsilon$, we have that

$$(47) \quad X_i(1) > \frac{M}{K} = \frac{M}{2\alpha Mc + 1} \quad \text{for all } i \in [-L/2, L/2]^d \cap \mathbb{Z}^d$$

and

$$(48) \quad \sup_{1 \le t \le 2} \{Y_j(t) : j \in [-NL, NL]^d \cap \mathbb{Z}^d\} < a'M',$$

where $N = b + 2$. Thus, throughout the time interval $[1, 2]$, a "safe environment condition" $G'_{\text{sec}}(i, 1) \in \mathcal{F}'(i, NL, 1, 1)$ holds at time $n = 1$ for all $i \in [-(3/2)L, (3/2)L]^d \cap \mathbb{Z}^d$ and hence the local $Y$-population is not "too big," that is, for all $j \in [-(3/2)L, (3/2)L]^d \cap \mathbb{Z}^d$, we have $\max_{t \in [1,2]} \sum_{l \in \mathbb{Z}^d} \gamma_{jl} Y_l(t) < 1$. Thus, the $\zeta$-process can safely invade the neighboring boxes, that is, conditional on (47) and given the above instance of the "safe environment condition," for each site

$$(49) \quad i \in \{\{L\mathbf{e_1} + [-L/2, L/2]^d\} \cup \{-L\mathbf{e_1} + [-L/2, L/2]^d\}\} \cap \mathbb{Z}^d,$$

where $\mathbf{e_1}$ denotes the first unit vector in $\mathbb{Z}^d$, the "infection event" $G_{\text{infec}}(i, j, 1)$ at $i$ has probability greater than $1 - \varepsilon$, by Lemma 4.3 (noting that by our choice of $L$, and the fact that $m_{ij}$ is a function of $\|i - j\|$ alone, each such site $i$ has at least one occupied neighbor $j \in [-L/2, L/2]$). Hence, after all of these prerequisites, the probability that simultaneously for all such sites $i$ taken from the set in (49), at time 1, the event

$G_{\text{infec}}(i, j, 1) \in \mathcal{F}(i, NL, 1, 1)$ holds for some $j \in [-L/2, L/2]^d \cap \mathbb{Z}^d$,

implying, under the above conditions that $X_i(2) > M/K = M/(2\alpha Mc + 1)$, is at least $1 - 2L^d \varepsilon$. We denote this event by $G_{\text{infec}}(L \uparrow L/2, 1, 1) \in \mathcal{F}(0, NL, 0, 2)$. Thus, we may define the $\mathcal{F}^*(NL, [0, 2])$-measurable "good event"

$$(50) \quad \begin{aligned} G_\zeta := \ & E'(0, (b+1/2)L, 0, 2)[a'] \cap G_{\text{nonrec}}(i, L/2, 0, 1) \\ & \cap E'(i, NL, 0, 2)[a'/2] \cap G_{\text{infec}}(i, L \uparrow L/2, 1, 1), \end{aligned}$$

which implies, given $\zeta = \zeta(0) \in H$, that $\zeta(1) \in (+1)H \cap (-1)H$ and observe

$$\mathbb{P}[G_\zeta] > 1 - (2(b+1/2)L)^d \varepsilon - L^d \varepsilon - (2NL)^d \varepsilon - 2L^d \varepsilon$$
$$> 1 - 4(2NL)^d \varepsilon = 1 - \theta,$$

which completes the comparison. $\square$

4.2.2. *Simultaneous comparison and proof of Theorem* 1.2 *and Corollary* 1.3. Assume that $\theta \leq 6^{-4(4(b+2)+1)^2}$. We may then choose $a, a'$ such that for all $M > \sum_j m_{ij}/\alpha$ and $M' > \sum_j m'_{ij}/\alpha'$, (45) holds for both populations $X$ and $Y$ with $\varepsilon = \frac{1}{4}\frac{\theta}{(2NL)^d}$. The point is that this can be done simultaneously, since the bounds for the control of the environment do not depend on the behavior of the competitor.

We may then pick $M, M'$ and simultaneously $\gamma, \gamma'$ such that Lemma 4.3 holds with $\varepsilon = \frac{1}{4}\frac{\theta}{(2NL)^d}$ for both the $X$- and the $Y$-population. Condition (45) is unaffected by this, since the bounds on the environment do not depend on $\{\gamma_{ij}\}$ and $\{\gamma'_{ij}\}$ and hold for all $M > \sum_j m_{ij}/\alpha$ and $M' > \sum_j m'_{ij}/\alpha'$. Assuming, then, that $\zeta(0)$ and $\eta(0) \in H$, observing that condition (44) on $a'$ ensures that $\kappa'_1 := \frac{M'}{2\alpha'M'c+1} < a'M' =: \kappa'_2$ (with a similar inequality for $\kappa_1, \kappa_2$) leads to the initial condition

$$(X(0), Y(0)) \in H(\kappa_1, \kappa_2; \kappa'_1, \kappa'_2; (b+1/2)L)$$

specified in Theorem 1.2, with

$$\kappa_1 = \frac{M}{2\alpha Mc+1}, \qquad \kappa_2 = aM,$$
$$\kappa'_1 = \frac{M'}{2\alpha'M'c+1}, \qquad \kappa'_2 = a'M'.$$

Hence, we can simultaneously construct the corresponding good events $G_\zeta$ and $G_\eta$ and infer from Theorem 3.2 that both the $X$- and $Y$-population survive, each with probability greater than $\frac{19}{20}$, which yields *persistence* of $\{X, Y\}$ with positive probability. Moreover, if we make the stronger assumption that the initial configurations of the $X$- and $Y$-populations satisfy

$$(X(0), Y(0)) \in H(\kappa_1, \kappa_2; \kappa'_1, \kappa'_2; \infty),$$

hence assuming $\zeta_i(0) = 1$, $\eta_i(0) = 1$, for all $i \in \mathbb{Z}^d$, then, according to Theorem 3.6, $\liminf_{n \to \infty} \mathbb{P}[\zeta_{2n}(0) = 1] \geq \frac{19}{20}$. The same result holds for $\eta$. Thus,

$$\liminf_{n \to \infty} \mathbb{P}[\zeta_{2n}(0) = 1, \eta_{2n}(0) = 1] \geq \tfrac{9}{10}.$$

By the definition of $\{\zeta_n, \eta_n\}_{n \geq 0}$ and our bounds from the last section applied in a similar fashion, this implies that there is a uniform positive lower bound on $\mathbb{P}[X_0(t) \geq M/(2K), Y_0(t) \geq M/(2K)\ \forall t \in [0, 4]$ given $\zeta_{2n}(0) = 1, \eta_{2n}(0) = 1]$ and the proof is completed.

4.2.3. *Proof of Theorem* 1.4. The proof of Theorem 1.4 again follows by comparison to suitable one-dimensional diffusions. This time, there is no need to control a potentially harmful "environment," making things much easier. For details, we again refer to Blath, Etheridge and Meredith [1].

# APPENDIX

The following classical comparison theorem and, in particular, its subsequent corollary, whose proof can be found in Blath, Etheridge and Meredith [1], are tailored for our purposes in Section 4. Note that the corollary allows a comparison, even if (51) only holds for intervals.

THEOREM A.1 (Ikeda and Watanabe [10]). *Let $(\Omega, \mathcal{F}, \{\mathcal{F}_t\}, \mathbb{P})$ be a filtered probability space and let $x_1(t,\omega), x_2(t,\omega)$ be two real $\{\mathcal{F}_t\}$-adapted processes. Let $B(t,\omega)$ be a one-dimensional $\{\mathcal{F}_t\}$-Brownian motion such that $B(0) = 0$ a.s. and let $\beta_1(t,\omega), \beta_2(t,\omega)$ be two real $\{\mathcal{F}_t\}$-adapted previsible drifts. Assume that with probability one,*

$$x_i(t) - x_i(0) = \int_0^t \sqrt{x_i(t)}\, dB(s) + \int_0^t \beta_i(s)\, ds, \qquad i = 1, 2,$$

*and that pathwise uniqueness of solutions holds for at least one of the equations. Moreover, assume that with probability one,*

$$x_1(0) \leq x_2(0), \qquad \beta_1(t) \leq b_1(t, x_1) \quad and \quad \beta_2(t) \geq b_2(t, x_2) \qquad \forall t \geq 0$$

*for two real continuous functions $b_1(t,x), b_2(t,x)$ on $[0,\infty) \times \mathbb{R}$ such that*

$$b_1(t,x) \leq b_2(t,x) \tag{51}$$

*for all $t \geq 0$ and $x \in \mathbb{R}$. Then, $x_1(t) \leq x_2(t)$ for every $t \geq 0$.*

COROLLARY A.2. *In the framework of Theorem A.1, assume that $x_1$ and $x_2$ are positive and nonexploding. Let $\delta > 0$.*

(a) *Suppose that condition (51) on $b_1, b_2$ is required only for all $x \in [\delta, \infty)$. Assume $x_1(0) \leq x_2(0)$. Define $\tau_\delta^{x_2} := \inf\{t \geq 0 : x_2(t) \leq \delta\}$. Then, with probability one, $x_1(t \wedge \tau_\delta^{x_2}) \leq x_2(t \wedge \tau_\delta^{x_2})$ for all $t \geq 0$.*
(b) *Suppose that condition (51) on $b_1, b_2$ is required only for all $x \in [0, \delta]$. Assume $x_1(0) \leq x_2(0)$. Define $\tau_\delta^{x_1} := \inf\{t \geq 0 : x_1(t) \geq \delta\}$. Then, with probability one, $x_1(t \wedge \tau_\delta^{x_1}) \leq x_2(t \wedge \tau_\delta^{x_1})$ for all $t \geq 0$.*

**Acknowledgments.** We wish to thank Matthias Birkner for helpful discussions. We are very grateful to two anonymous referees for very careful reading of the manuscripts and for many improvements to the presentation.

## REFERENCES

[1] BLATH, J., MEREDITH, M. E. and ETHERIDGE, A. M. (2007). Coexistence in locally regulated competing populations and survival of BARW: Full technical details and additional remarks. Preprint, TUB.
[2] BOLKER, B. M. and PACALA, S. W. (1999). Spatial moment equations for plant competition: Understanding spatial strategies and the advantages of short dispersal. *American Naturalist* **153** 575–602.

[3] CARDY, J. L. and TÄUBER, U. C. (1996). Theory of branching and annihilating random walks. *Phys. Rev. Lett.* **77** 4780–4783.
[4] CARDY, J. L. and TÄUBER, U. C. (1998). Field theory of branching and annihilating random walks. *J. Statist. Phys.* **90** 1–56. MR1611125
[5] COX, J. T. and PERKINS, E. A. (2005). Rescaled Lotka–Volterra models converge to super-Brownian motion. *Ann. Probab.* **33** 904–947. MR2135308
[6] COX, J. T. and PERKINS, E. A. (2006). Survival and coexistence in stochastic spatial Lotka–Volterra models. Preprint.
[7] DIECKMANN, U. and LAW, R. (2000). Relaxation projections and the method of moments. In *The Geometry of Ecological Interactions*: *Simplifying Spatial Complexity* (U. Dieckmann, R. Law and J. A. J. Metz, eds.) 412–455. Cambridge Univ. Press.
[8] DURRETT, R. (1995). Ten lectures on particle systems. *Ecole d'Été de Probabilités de Saint Flour XXIII 1993. Lecture Notes in Math* **1608** 97–201. Springer, Berlin. MR1383122
[9] ETHERIDGE, A. M. (2004). Survival and extinction in a locally regulated population. *Ann. Appl. Probab.* **14** 188–214. MR2023020
[10] IKEDA, N. and WATANABE, S. (1981). *Stochastic Differential Equations and Diffusion Processes.* North-Holland and Kodensha, Tokyo. MR0637061
[11] LAW, R., MURRELL, D. J. and DIECKMANN, U. (2003). Population growth in space and time: Spatial logistic equations. *Ecology* **84** 252–262.
[12] MURRELL, D. J. and LAW, R. (2003). Heteromyopia and the spatial coexistence of similar competitors. *Ecology Letters* **6** 48–59.
[13] NEUHAUSER, C. and PACALA, S. W. (1999). An explicitly spatial version of the Lotka–Volterra model with interspecific competition. *Ann. Appl. Probab.* **9** 1226–1259. MR1728561
[14] PACALA, S. W. and LEVIN, S. A. (1997). Biologically generated spatial patterns and the coexistence of competing species. In *Spatial Ecology*: *The Role of Space in Population Dynamics and Interspecific Interactions* (D. Tilman and P. Kareiva, eds.) 204–232. Princeton Univ. Press.
[15] SHIGA, T. (1980). An interacting system in population genetics. *J. Math. Kyoto Univ.* **20** 213–242. MR0582165
[16] SHIGA, T. (1982). Continuous time multi-allelic stepping stone models in population genetics. *J. Math. Kyoto Univ.* **22** 1–40. MR0648554
[17] SHIGA, T. and SHIMIZU, A. (1980). Infinite dimensional stochastic differential equations and their applications. *J. Math. Kyoto Univ.* **20** 395–416. MR0591802
[18] TÄUBER, U. C. (2003). Scale invariance and dynamic phase transitions in diffusion-limited reactions. In *Advances in Solid State Physics* (B. Kramer, ed.) 659–675. Springer, Berlin.

J. BLATH
INSTITUT FÜR MATHEMATIK
TECHNISCHE UNIVERSITÄT BERLIN
STRASSE DES 17. JUNI 136
10623 BERLIN
GERMANY
E-MAIL: blath@math.tu-berlin.de

E. ETHERIDGE
M. MEREDITH
DEPARTMENT OF STATISTICS
UNIVERSITY OF OXFORD
1 SOUTH PARKS ROAD
OXFORD OX1 3TG
UNITED KINGDOM
E-MAIL: alison.etheridge@magd.ox.ac.uk
  mark.meredith@magd.ox.ac.uk